\documentclass[11pt]{article}
\usepackage{amsfonts}
\usepackage{lineno,hyperref}
\pdfoutput=1
\usepackage{bbm}
\usepackage{amsmath}
\usepackage{mathrsfs}
\usepackage{amssymb}
\usepackage{amsmath,color}
\usepackage{amsthm}
\usepackage{amstext}
\usepackage{amsopn}
\usepackage{texdraw}
\usepackage{graphicx}
\usepackage{multirow}
\usepackage{epsfig,subfig,color}
\usepackage{lscape}
\usepackage{rotating}
\usepackage{float}
\usepackage{cite}
\usepackage{srcltx}
\usepackage[square, comma, sort&compress, numbers]{natbib}
\usepackage{appendix}
\DeclareGraphicsExtensions{.pdf,.eps,.png,.jpg,.mps}
\newtheorem{theorem}{Theorem}[section]

\newtheorem{proposition}[theorem]{Proposition}

\theoremstyle{definition}

\theoremstyle{remark}
\newtheorem{remark}[theorem]{Remark}
\numberwithin{equation}{section}

\graphicspath{{figequilibrium/}}
\date{}

\begin{document}
\title{ \bf\large{The spatially inhomogeneous Hopf bifurcation induced by memory delay  in a memory-based diffusion system}\footnote{Partially supported by the Natural Science Foundation of  China (Nos.11971143, 12071105), and  Natural Science Foundation of Zhejiang Province of China  (No.LY19A010010).}}
\author{Yongli Song\textsuperscript{1}, \  \ Yahong Peng\textsuperscript{2}\footnote{Corresponding author, Email: pengyahong@dhu.edu.cn },  \  \ Tonghua Zhang\textsuperscript{3}  \\
{\small \textsuperscript{1} Department of Mathematics, Hangzhou Normal University, Hangzhou 311121, China\hfill {\ }}\\
{\small \textsuperscript{2} Department of Mathematics, Donghua University, Shanghai,
201620, China\hfill{\ }}\\
{\small \textsuperscript{3} Department of Mathematics, Swinburne University of Technology, Hawthorn, VIC3122,  Australia\hfill {\ }}\\
}
\maketitle
\begin{abstract}

The memory-based diffusion systems have wide applications in practice. Hopf bifurcations are observed from such systems. To meet the demand for computing the normal forms of the Hopf bifurcations of such systems, we develop an effective new algorithm where the memory delay is treated as the perturbation parameter. To illustrate the effectiveness of the algorithm, we consider a diffusive predator-prey system with memory-based diffusion and Holling type-II functional response. By employing this newly developed procedure, we   investigate the direction and stability of the delay-induced mode-1 and mode-2 Hopf bifurcations. Numerical simulations confirm our theoretical findings, that is  the existence of stable spatially inhomogeneous periodic solutions with mode-1 and mode-2 spatial patterns,  and the transition  from the  unstable mode-$2$ spatially inhomogeneous periodic solution to the stable mode-$1$ spatially inhomogeneous periodic solution.
\end{abstract}



\section{Introduction}\label{sec1}

The study of delay-induced Hopf bifurcations is a very active research area in recent years, and has attracted many researchers' attention due to its importance. In the study of nonlinear dynamical systems that may exist rich dynamics, such as Hopf bifurcation, the normal form plays an import role as it can be used to determine the direction and stability of the Hopf bifurcation. Thus, many researchers have dedicated themselves to the development of ways calculating normal forms of the nonlinear dynamical systems. For example,  Hassard, Kazarinoff and Wan \cite{Hassard-book} and Faria and Magalhaes \cite{Faria-JDE1995} proposed an algorithm for computing the normal form of the Hopf bifrucations for delay differential systems (without diffusion), respectively; Faria \cite{Faria-00TAMS} developed an algorithm for reaction-diffusion systems with time delays, but which only appear in the reaction terms; most recently Song and his colleagues \cite{Wu-Song-CNSNS2020,Song-Shi-2021} developed some algorithms for reaction-diffusion systems with both time delay and nonlocal spatial average. It is now known that in the classical delay reaction-diffusion systems, the first delay-induced Hopf bifurcation is often homogeneous and the bifurcating periodic oscillation is spatially homogeneous \cite{Zuo-Wei-NARWA2011,Zuo-Wei-NAMC2014,Zhang-Yan-IJBC2014,Tang-Song-AMC2015,Yang-Liu-Zhang-CNSNS2017,Tang-Jiang-JAAC2017,Yang-Zhang-IJBC2018,Song-Jiang-Yuan-JAAC2019,Jiang-Wang-Cao-JDDE2019,Shen-Wei-IJBC2019,Manna-Malay-MBE2019,Chen-WuLC-CNSNS2019,Chen-Zhang-Zhou-MCS2020,Dai-Liu-Wei-MBE2020,Meng-Wang-JBD2020,Li-Mei-Cao-IJBC2020,Lin-Xu-Li-CNSNS2020,Jiang-Zhang-She-IJB2020}.  And the aforementioned algorithms were developed for such scenarios.

However, in applications, the delay reaction-diffusion systems may not be homogeneous in space. For example, very recently,  to better understand  how memory influences the animal movement, Shi {\it et al.}  \cite{Shi-WangYan-JDDE2020}  proposed a modified Fick's law where in addition to the negative gradient of the density distribution function at the present time, there is a directed movement toward the negative gradient of the density distribution function at past time, and then based on the modified Fick's law, they formulated a single species diffusive model with spatial memory. Since then, there have recently been an increasing activities and interests on the study of the dynamics of the reaction-diffusion equations with spatial memory \cite{Shi-Wang-Wang-Nonlinearity2019,Shi-WangYan-JDDE2020,An-Wang-Wang-DCDS2020,Oliveira-Kar-Ber-MB2020,Shi-Shi-Wang-JMB2021,Song-Shi-Wang-2021}. For this kind of equations,  the phenomenon that the first delay-induced Hopf bifurcation is inhomogeneous has been observed by \cite{Shi-Wang-Wang-Nonlinearity2019,Song-Wu-Wang-JDE2019}.
 Based on the assumption that  the prey has no memory or cognition,  we introduced the spatial memory to the diffusive predator-prey model and proposed the  following  model with random and memory-based diffusions subject to Neumann boundary conditions on one-dimensional spatial domain $(0, \ell \pi)$ with $\ell\in \mathbb{R}^+$ is
\begin{equation}
\label{PPMODEL}
  \begin{cases}\frac{\partial u(x,t)}{\partial t}=d_{11}u_{xx}(x,t)+f(u(x,t),v(x,t)), & 0<x<\ell\pi, t>0, \vspace{0.1cm}\\
\frac{\partial v(x,t)}{\partial t}=d_{22}v_{xx}(x,t)-d_{21}(v(x,t)u_{x}(x,t-\tau))_x & \\
\hspace{0.6in} +g(u(x,t),v(x,t)), & 0<x<\ell\pi, t>0,\vspace{0.1cm}\\
u_x(0,t)=u_x(\ell\pi,t)=v_x(0,t)=v_x(\ell\pi,t)=0, &t\geq0,
  \end{cases}
 \end{equation}
where $d_{11}$ and $d_{22}$ are the Fickian diffusion coefficients, $d_{12}$ and $d_{21}$ are the memory-based diffusion coefficients, the time delay, denoted by $\tau>0$ represents the averaged memory period of predator. It has been shown that delay-induced first Hopf bifurcations are often spatially inhomogeneous \cite{Song-Shi-Wang-2021}.

Compared with the classical reaction-diffusion system with delay,  system \eqref{PPMODEL} has two main characteristics: $(i)$ the delay appears in the diffusion term; $(ii)$ the diffusion term is not linear.   Therefore, algorithms,  in the existing literatures such as \cite{Hassard-book,Faria-00TAMS,Wu-Song-CNSNS2020,Song-Shi-2021},  that were developed for analysing the Hopf bifurcation of the reaction-diffusion systems, can not be applied to the system \eqref{PPMODEL}. To bridge the gap, in this paper, we aim to develop an algorithm that can calculate the normal form for analysing the Hopf bifurcation which may occur in a spatially inhomogeneous scenario and be induced by the memory delay in \eqref{PPMODEL}. And then we apply this algorithm to a diffusive predator-prey system with memory-based diffusion and Holling type-II functional response. More precisely, we organize the rest of our paper as follows.
 In Section \ref{sec2}, we derive the algorithm of calculating normal form of spatially Hopf bifurcation  induced by memory delay for system \eqref{PPMODEL}.   In Section \ref{sec3}, we study  the direction and stability of delay-induced mode-1 and mode-2 Hopf bifurcations in the diffusive predator-prey system with memory-based diffusion and Holling type-II functional response and the numerical simulations are used to illustrate the theoretical results. Finally,  we give a brief discussion and conclude our study in Section \ref{sec4}. For the convenience of discussion, throughout the paper, we let $\mathbb{N}$ represent the set of natural numbers, and $\mathbb{N}_0=\mathbb{N}\cup \{0\}$ represent the set of all nonnegative integers.

\section{Calculation of the normal form of the Hopf bifurcation}\label{sec2}

\begin{subsection}
   {Characteristic equation at the equilibrium and some basic assumptions}
 \end{subsection}
 Without loss of generality, assume that $E_*=(u_*, v_*)$ is a constant coexistence (positive) equilibrium of system \eqref{PPMODEL}. Then straightforward calculation gives the linearization of \eqref{PPMODEL} at $(u_*, v_*)$:
\begin{equation}
\label{LPPMODEL}
\left(\begin{array}{c}
u_t(x,t)\\
v_t(x,t)
\end{array}\right)=D_1\left(\begin{array}{c}
u_{xx}(x,t)\\
v_{xx}(x,t)
\end{array}\right)+D_2\left(\begin{array}{c}
u_{xx}(x,t-\tau)\\
v_{xx}(x,t-\tau)
\end{array}\right)+A\left(\begin{array}{c}
u(x,t)\\
v(x,t)
\end{array}\right),
 \end{equation}
 where
  \begin{equation}
 \label{DJA}
D_1=\left(\begin{array}{cc}
d_{11} & 0\\
0 & d_{22}
\end{array}\right),~D_2=\left(\begin{array}{cc}
0 & 0\\
-d_{21}v_* & 0
\end{array}\right),  ~A=\left(
                          \begin{array}{cc}
                            a_{11} & a_{12} \\
                            a_{21} & a_{22} \\
                          \end{array}
                        \right),
\end{equation}
and
 $$
 a_{11}=\frac{\partial f(u_*, v_*)}{\partial u },  ~a_{12}=\frac{\partial f(u_*, v_*)}{\partial v }, ~a_{21}=\frac{\partial g(u_*, v_*)}{\partial u },  ~a_{22}=\frac{\partial g(u_*, v_*)}{\partial v }.
 $$
Thus, the characteristic equation of \eqref{LPPMODEL} is
\begin{equation}
\label{CE0}
\prod\limits_{n\in \mathbb{N}_0} \Gamma_n(\lambda)=0,
\end{equation}
where $\Gamma_n(\lambda)=\mathrm{det}\left(\mathcal{M}_n(\lambda)\right)$. Notice
\begin{equation}
\label{CEN}
\mathcal{M}_n(\lambda)=\lambda I_2+(n/\ell)^2D_1+(n/\ell)^2e^{-\lambda\tau}D_2-A,
\end{equation}
 with $I_2$ is the identity matrix of size 2.
Then we have
\begin{equation}
\label{CE}
\Gamma_n(\lambda)=\mathrm{det}\left(\mathcal{M}_n(\lambda)\right)=\lambda ^2-T_n\lambda+\widetilde{J}_n(\tau)=0,
\end{equation}
where
\begin{equation}
\label{TK}
\begin{split}
T_n&=Tr(A)-Tr(D_1)(n/\ell)^2,\\
\widetilde{J}_n(\tau)&=d_{11}d_{22} (n/\ell)^4
-\left(d_{11}a_{22}+d_{22}a_{11}+d_{21}v_*a_{12}e^{-\lambda\tau}\right) (n/\ell)^2+Det(A),
\end{split}\end{equation}
with $Tr(A)=a_{11}+a_{22}$,  $Tr(D_1)=d_{11}+d_{22}$ and  $Det(A)=a_{11}a_{22}-a_{12}a_{21}$.

For the Hopf bifurcation, we assume that  at $\tau=\tau_c$, Eq.\eqref{CE} has a pair of purely imaginary  roots $\pm i\omega_{n_c}, \omega_{n_c}>0$  for $n=n_c\in\mathbb{N}$  and all other eigenvalues  have negative real part.  Let
$\lambda(\tau)=\alpha(\tau)+i\beta(\tau)$ be a pair of roots of
Eq.\eqref{CE} near $\tau=\tau_c$ satisfying
$\alpha\left(\tau_c\right)=0$ and $\beta\left(\tau_c\right)=\omega_{n_c}$. Then, it is easy to verify  the transversality condition
 \[
 \left.\frac{d \mathrm{Re}(\lambda (\tau))}{d
\tau}\right|_{\tau=\tau_c}\neq 0.
\]
\begin{subsection}
   {Algorithm for computing the normal form}
 \end{subsection}

We first introduce some conventional notations: define the  two-dimensional real-valued Sobolev
space
\begin{equation*}
\mathscr{X}=\left\{U=\left(U_1, U_2\right)^T\in \left(W^{2,2}(0, \ell\pi)\right)^2,\frac{\partial U_1}{\partial
x}=0, \frac{\partial U_2}{\partial
x}=0,~ x=0, \ell\pi\right\},
\end{equation*}
with the inner product defined by
\begin{equation*}
\left[U, V \right]=\int_0^{\ell\pi}U^T V
dx, ~\mbox{for}~  U, V \in \mathscr{X},
\end{equation*}
and let  $\mathscr{C}:=C\left([-1,0]; \mathscr{X}\right)$ be the Banach space of continuous mappings from $[-1, 0]$ to $\mathscr{X}$ .

Then, we take a small perturbation of $\tau_c$ by setting $\tau=\tau_c+\mu, |\mu|\ll1$ such that $\mu=0$ corresponds to the Hopf bifurcation value for Eq.\eqref{PPMODEL}. Here, $D_j, j=1, 2$,  are given by \eqref{DJA}.

Now shift $E_*$ to the origin by setting
$$U(x,t)=\left(U_1(x,t),  U_2(x,t)\right)^T=(u(x,t), v(x,t))^T-(u_*, v_*)^T,$$
normalize the delay by rescaling the time variable $t\rightarrow t/\tau$, and rewrite $U(t)$ for $U(x, t)$ and $U_t\in \mathscr{C}$ for $U_t(\theta)=U(x, t+\theta), -1\leq \theta\leq 0$. Then, the system \eqref{PPMODEL} becomes the compact form
\begin{equation}
\label{NFEQ}
\frac{dU(t)}{dt}=d(\mu)\Delta (U_t)+L(\mu)(U_t)+F(U_t, \mu),
\end{equation}
where   for $\varphi=\left(\varphi^{(1)},\varphi^{(2)}\right)^T\in \mathscr{C}$,  $d(\mu)\Delta, L(\mu): \mathscr{C}\rightarrow \mathscr{X}$, $F: \mathscr{C}\times \mathbb{R}^2\rightarrow \mathscr{X}$  are given, respectively, by
$$d(\mu)\Delta (\varphi)=d_0\Delta (\varphi)+F^d (\varphi, \mu), ~L(\mu) (\varphi)=(\tau_c+\mu) A\varphi(0),$$
 and
\begin{equation}
\label{FPM}
  F(\varphi, \mu) = \left(\tau_c+\mu\right) \left(\begin{array}{c}
f\left(\varphi^{(1)}(0)+u_*, \varphi^{(2)}0)+v_*\right)\vspace{0.2cm}\\
g\left(\varphi^{(1)}(0)+u_*, \varphi^{(2)}(0)+v_*\right)\end{array}\right)-L(\mu)(\varphi),
\end{equation}
where
 $$
 d_0\Delta (\varphi)=\tau_cD_1\varphi_{xx}(0)+\tau_cD_2^c\varphi_{xx}(-1),
$$
\begin{equation}
\label{FD}
\begin{array}{lll}
F^d (\varphi, \mu)
&=&
-d_{21}\left(\tau_c+\mu\right)\left(\begin{array}{c}
0\vspace{0.2cm}\\
\varphi^{(1)}_{x}(-1)\varphi^{(2)}_{x}(0)+\varphi^{(1)}_{xx}(-1)\varphi^{(2)}(0)
\end{array}\right)\vspace{0.3cm}\\
&&
+\mu \left(\begin{array}{c}
d_{11}\varphi^{(1)}_{xx}(0)\vspace{0.2cm}\\
-d_{21} v_*\varphi^{(1)}_{xx}(-1) +d_{22}\varphi^{(2)}_{xx}(0)\end{array}\right).
\end{array}
 \end{equation}
  In what follows, we assume that  $F(\varphi, \mu)$ is $C^k, k\geq 3,$ smooth with respect to $\varphi$ and $\mu$.
  Noticing that $\mu$ is the perturbation parameter and treated as a variable in the calculation of normal forms,  we denote
  $L_0(\varphi)=\tau_cA\varphi(0)$
 and rewrite  \eqref{NFEQ}   as  the following
  \begin{equation}\label{NFEQ1}
\frac{dU(t)}{dt}=d_0\Delta (U_t)+L_0(U_t)+\widetilde{F}(U_{t},  \mu),
\end{equation}
by separating the linear terms from the nonlinear terms,
\begin{equation}
\label{NLTERM}
\widetilde{F}(\varphi, \mu)=\mu A\varphi(0)+F(\varphi, \mu)+F^d(\varphi, \mu).
\end{equation}

Then the characteristic equation for the linearized system of \eqref{NFEQ1}
 \begin{equation}\label{LNFEQ1}
\frac{dU(t)}{dt}=d_0\Delta (U_t)+L_0(U_t)
\end{equation}
is
\begin{equation}
\label{NFCE}
\prod\limits_{n\in \mathbb{N}_0} \widetilde{\Gamma}_n(\lambda)=0,
\end{equation}
where  $\widetilde{\Gamma}_n(\lambda)=\mathrm{det}\left(\widetilde{\mathcal{M}}_n(\lambda)\right)$ with
\begin{equation}
\label{NFMN}
\widetilde{\mathcal{M}}_n(\lambda)=\lambda I_2+\tau_c(n/\ell)^2D_1+\tau_c(n/\ell)^2e^{-\lambda}D_2-\tau_cA.
\end{equation}

Comparing   \eqref{NFMN} with \eqref{CEN}, we know that
Eq.\eqref{NFCE}  has a pair  of purely imaginary  roots
$\pm i\omega_c$   for $n=n_c\in\mathbb{N}$,  and  all other eigenvalues  have negative real parts, where
$\omega_c=\tau_c\omega_{n_c}$.

 It is well known that the eigenvalue problem
 $$
-\gamma''=\mu \gamma,~~x\in(0, \ell\pi);~~\gamma'(0)=\gamma'(\ell\pi)=0
$$
\ has eigenvalues $\mu_{n}=(n/\ell)^{2}, n\in \mathbb{N}_0 $,
with the corresponding normalized eigenfunctions
$$\gamma_n(x)=\frac{\cos\left(\frac{nx}{\ell}\right)}{\|\cos\left(\frac{nx}{\ell}\right)\|_{2,2}}=\left\{
\begin{array}{ll}
  \frac{1}{\sqrt{\ell \pi}}, & ~\mathrm{for}~n=0,\vspace{0.2cm} \\
  \frac{\sqrt{2}}{\sqrt{\ell \pi}}\cos\left(\frac{nx}{\ell}\right), & ~\mathrm{for}~n\neq 0.
\end{array}
\right.$$

Let  $\beta_n^{(j)}= \gamma_n(x)e_j, j=1, 2$,
where $e_j$ are the unit coordinate vectors of $\mathbb{R}^2$. Then the eigenfunctions $\left\{\beta_n^{(j)}\right\}_{n=1}^{\infty}$  form an orthonormal basis for $\mathscr{X}$.

Let $\mathscr{B}_n=\mbox{span}\left\{\left[v(\cdot), \beta_n^{(j)}\right]\beta_n^{(j)}|~ v\in \mathscr{C}, j=1, 2\right\}$.  Then it is easy to verify that
$$L_0(\mathscr{B}_n)\subset \mbox{span}\left\{\beta_n^{(1)}, \beta_n^{(2)}\right\}, n\in \mathbb{N}_0.$$
Assume that $z_t(\theta)\in C=C\left([-1,0], \mathbb{R}^2\right)$ and
 $$z_t^T(\theta) \left(\begin{array}{c}
                             \beta_n^{(1)} \vspace{0.1cm}\\
                             \beta_n^{(2)}
                           \end{array}
\right)\in \mathscr{B}_n.$$
Then, on $\mathscr{B}_n$, the linearized equation  \eqref{LNFEQ1} is equivalent to the following  functional differential equation (FDE) in $C$:
\begin{equation}
\label{FDE}
\dot{z}(t)=L_0^d\left(z_t(\theta) \right)+ L_0(z_t(\theta)),
\end{equation}
where
$$
L_0^d\left(z_t(\theta) \right) =\tau_c\left(\begin{array}{cc}
                  -d_{11} (n/\ell)^2 & 0 \\
                   0 &   -d_{22}(n/\ell)^2
                 \end{array}
\right) z_t(0)+\tau_c\left(\begin{array}{cc}
                  0 & 0 \\
                d_{21}v_*(n/\ell)^2 &  0
                 \end{array}
\right) z_t(-1).
$$
The characteristic equation of  linear system  \eqref{FDE} is the same as the one  given in  \eqref{NFCE}.

Define $\eta_n(\theta)\in BV([-1, 0], \mathbb{R}^2)$
such that
$$
\int^0_{-1}d\eta_n(\theta) \varphi(\theta)=L_0^d(\varphi(\theta)) + L_0(\varphi(\theta)), ~~\varphi\in C,
 $$
and use the adjoint bilinear form on $C^*\times C$, $C^*=C([0, 1], \mathbb{R}^{2*})$, where $\mathbb{R}^{2*}$ is a 2-dimensional
space of row vectors,   as follows
$$
\langle\psi(s),
\varphi(\theta)\rangle_n =\psi(0)\varphi(0)-\int^0_{-1}\int^{\theta}_0
\psi(\xi-\theta)d\eta_n (\theta)\varphi(\xi)d\xi,\ \ \
\textmd{for}\ \psi\in C^*, \varphi\in C.
$$

Let
$\Lambda=\{i\omega_c, -i\omega_c\}$. Denote the generalized eigenspace of \eqref{FDE} associated with $\Lambda$  by $P$ and  the corresponding adjoint space by
$P^*$. Then, by the adjoint theory of functional differential equation \cite{HALE-1977BOOK}, $C$ can be decomposed as
$C=P\oplus Q$, where $Q=\{\varphi\in C:
\langle\psi,\varphi\rangle=0,\forall\psi\in P^*\}$. Choose the
bases $\Phi(\theta)$ and  $\Psi(s)$ of $P$ and $P^*$, respectively, as follows
$$\Phi(\theta)=\left(\phi(\theta), \overline{\phi}(\theta)\right),
 ~\Psi(s)=\mbox{col}\left(\psi^T(s), \overline{\psi }^T(s)\right),$$
such that $\langle\Psi, \Phi\rangle_{n_c}=I_2,$
where
$$
\phi(\theta)=\left(
\begin{array}{c}
\phi_1(\theta)\vspace{0.2cm}\\
 \phi_2(\theta)
 \end{array}
\right)=\phi e^{i\omega_{c}\theta},~\psi(s)=\left(
\begin{array}{c}
\psi_1(\theta)\vspace{0.2cm}\\
 \psi_2(\theta)
 \end{array}
\right)=
\psi e^{-i\omega_{c}s},
$$
and
\begin{equation*}
\phi=\left(\begin{array}{c} 1 \vspace{0.2cm}\\
\frac{i\omega_{n_c}+(n_c/\ell)^2d_{11}-a_{11}}{a_{12}}
\end{array}\right),
~~
\psi=\eta \left(\begin{array}{c} 1 \vspace{0.2cm}\\
\frac{a_{12}}{i\omega_{n_c}+(n_c/\ell)^2d_{22}-a_{22}}
\end{array}\right),
\end{equation*}
with
$$\eta =\frac{i\omega_{n_c}+(n_c/\ell)^2d_{22}-a_{22}}{i\omega_{n_c}+(n_c/\ell)^2d_{11}-a_{11}+i\omega_{n_c}+(n_c/\ell)^2d_{22}-a_{22}+\tau_ca_{12}d_{21}v_*(n_c/\ell)^2e^{-i\omega_c}}.$$

Using the decomposition $C=P\oplus Q$,  the phase space $\mathscr{C}$ for \eqref{NFEQ}  can be decomposed as
$$ \mathscr{C}=\mathcal{P} \oplus \mathcal{Q}, ~\mathcal{P}=\mbox{Im} \pi, ~\mathcal{Q}=\mbox{Ker} \pi, $$
 where $\pi:  \mathscr{C}\rightarrow \mathcal{P}$ is the  projection operator defined by
$$
\begin{array}{lll}
  \pi (\phi)
&=&
\Phi (\theta) \left \langle  \Psi (0), \left(\begin{array}{c}\left[\phi(\cdot), \beta_{n_c}^{(1)}\right]\vspace{0.2cm}\\
\left[\phi(\cdot),   \beta_{n_c}^{(2)}\right]
\end{array}\right)  \right\rangle  \gamma_{n_c}(x).
\end{array}
$$

Following \cite{Faria-00TAMS}  and \cite{Song-ZP-16CNSNS}, we
define $\mathscr{C}_0^1=\left\{ \phi\in \mathscr{C}: \dot{\phi}\in \mathscr{C}, \phi (0)\in \mbox{dom}(d\Delta)\right\}$ and let
$$z_x=\left(z_1(t)\gamma_{n_c}(x), z_2(t)\gamma_{n_c}(x)\right)^T, \Phi(\theta)=\left(\phi_{n_c}(\theta), \overline{\phi}_{n_c}(\theta)\right).
$$
For  $\varphi (\theta)\in \mathscr{C}_0^1$,   we have the following decomposition
$$
\varphi(\theta)=\Phi(\theta) z_x+w,  ~~w=(w^{(1)}, w^{(2)})^T\in \mathscr{C}_0^1\cap \mbox{Ker}\pi:=\mathscr{Q}^1.
$$
Following the notations in \cite{Faria-00TAMS}, we define
$$X_0(\theta)=\left\{
\begin{array}{ll}
0, & -1\leq \theta <0,\vspace{0.2cm}\\
1, & \theta=0.
\end{array}\right.
$$
and then
\begin{equation}
\label{PIX0}
  \pi\left(X_0(\theta)\widetilde{F}_2\left(\Phi(\theta) z_x,0\right)\right)\vspace{0.2cm}\\
=
\Phi (\theta)\Psi (0)\left(\begin{array}{c}\left[\widetilde{F}_2\left(\Phi(\theta) z_x,0\right),\beta_{n_c}^{(1)}\right]\vspace{0.2cm}\\
\left[\widetilde{F}_2\left(\Phi(\theta) z_x,0\right),\beta_{n_c}^{(2)}\right]
\end{array}\right)\gamma_{n_c}(x).
\end{equation}

Let $z=(z_1(t),  z_2(t))^T.$ Then system \eqref{NFEQ1} can be decomposed as a system of abstract
ODEs on $\mathbb{R}^4\times \mbox{Ker}\pi$:
\begin{equation}
\label{AODE}
\begin{cases}
\dot{z}=Bz+\Psi(0)\left(\begin{array}{c}\left[\widetilde{F}\left(\Phi(\theta) z_x+w,\mu\right),\beta_{n_c}^{(1)}\right]\vspace{0.1cm}\\
\left[\widetilde{F}\left(\Phi(\theta) z_x+w,\mu\right),\beta_{n_c}^{(2)}\right]
\end{array}\right),\vspace{0.1cm}\\
\dot{w}=A_{\mathcal{Q}^1}w+(I-\pi)X_0(\theta)\widetilde{F}\left(\Phi(\theta) z_x+w,\mu\right),
\end{cases}
\end{equation}
where
$B=\mbox{diag}\left\{i\omega_c, -i\omega_c  \right\},$
$A_{\mathcal{Q}^1}: \mathcal{Q}^1\rightarrow \mbox{Ker}\pi$ is defined by
$$A_{\mathscr{Q}^1} w=\dot{w}+ X_0(\theta) \left(L_0(w)+L_0^d(w)-\dot{w}(0) \right).$$

Consider the formal Taylor expansion
$$\widetilde{F}(\varphi, \mu)=\sum\limits_{j\geq 2}\frac{1}{j!}\widetilde{F}_j(\varphi, \mu),  ~~ F(\varphi, \mu)=\sum\limits_{j\geq 2}\frac{1}{j!} F_j(\varphi, \mu),  ~~ F^d(\varphi, \mu)=\sum\limits_{j\geq 2}\frac{1}{j!} F_j^d(\varphi, \mu).$$
From \eqref{NLTERM}, we have
\begin{equation}
\label{TLDF2MU}
\widetilde{F}_2( \varphi, \mu)=2\mu A\varphi(0)+
F_2(\varphi, \mu)+F_2^d(\varphi, \mu)
\end{equation}
and
\begin{equation}
\label{TLDF3MU}
\widetilde{F}_j( \varphi, \mu)=
F_j(\varphi, \mu) +F_j^d(\varphi, \mu), j=3, 4, \cdots.
\end{equation}

Then \eqref{AODE} is written as
$$
\left\{
\begin{array}{l}
  \dot{z}=Bz+\sum\limits_{j\geq 2} \frac{1}{j!}f^1_j(z, w, \mu),\vspace{0.2cm}\\
 \dot{w}=A_{\mathscr{Q}^1}w+ \sum\limits_{j\geq 2} \frac{1}{j!}f^2_j(z, w, \mu),
\end{array}
\right.
$$
where
\begin{equation}
\label{FJ1}
 f^1_j(z, w, \mu)=\Psi(0)\left(\begin{array}{c}\left[\widetilde{F}_j\left(\Phi(\theta) z_x+w,\mu\right),\beta_{n_c}^{(1)}\right]\vspace{0.1cm}\\
\left[\widetilde{F}_j\left(\Phi(\theta) z_x+w,\mu\right),\beta_{n_c}^{(2)}\right]
\end{array}\right),
\end{equation}
\begin{equation}
\label{FJ2}
 f^2_j(z, w, \mu)=(I-\pi)X_0(\theta) \widetilde{F}_j\left(\Phi(\theta) z_x+w, \mu\right).
\end{equation}

In terms of the normal form theory  of partial functional differential equations \cite{Faria-00TAMS}, after a recursive transformation of variables of the form
\begin{equation}
\label{TOV}
(z, w)=(\widetilde{z}, \widetilde{w})+\frac{1}{j!}\left(U_j^1(\widetilde{z}, \mu), U_j^2(\widetilde{z}, \mu)(\theta)\right), j\geq 2,
\end{equation}
where $z, \widetilde{z}\in \mathbb{R}^2, w, \widetilde{w}\in \mathscr{Q}^1$ and $U_j^1:\mathbb{R}^3\rightarrow  \mathbb{R}^2,U_j^2:\mathbb{R}^3\rightarrow  \mathscr{Q}^1$ are homogeneous polynomials of degree $j$ in $\widetilde{z}$ and $\mu,$
 the flow on the local center manifold for \eqref{NFEQ1} can be written as
\begin{equation}
\label{NF}
  \dot{z}=Bz+\sum\limits_{j\geq2}\frac{1}{j!}g_j^1(z,0,\mu),
\end{equation}
which is the normal form as in the usual sense for ODEs.

Following   \cite{Faria-00TAMS} and  \cite{Faria-JDE1995},   we have
$$
g_2^1(z, 0,\mu)=\mbox{Proj}_{\mbox{Ker}(M_2^1)}f_2^1(z,0,\mu),
$$
and
\begin{equation}
\label{G310}
g_3^1(z,0,\mu)=\mbox{Proj}_{\mbox{Ker}(M_3^1)}\widetilde{f}_3^1(z,0,\mu)=\mbox{Proj}_{S}\widetilde{f}_3^1(z,0,0)+O(\mu^2|z|),
\end{equation}
where  $\widetilde{f}_3^1(z,0,\mu)$ is vector, where elements are polynomials of degree 3 in $(z, \mu)$ obtained from \eqref{TOV} after performing the change of variables, and is determined by  \eqref{FTLD31},
\begin{equation*}
\mbox{Ker}\left(M_2^1\right)=\mbox{Span}\left\{
\left(\begin{array}{l}
 \mu z_1\\
0
\end{array}
\right),
\left(\begin{array}{l}
0\\
\mu z_2
\end{array}
\right)\right\},
 \end{equation*}
\begin{equation*}
\mbox{Ker}\left(M_3^1\right)=\mbox{Span}\left\{
\left(\begin{array}{l}
z_1^2z_2\\
0
\end{array}
\right),
 \left(\begin{array}{l}
\mu^2 z_1\\
0
\end{array}
\right),
\left(\begin{array}{l}
0\\
z_1z_2^2
\end{array}
\right),
 \left(\begin{array}{l}
0\\
\mu^2 z_2
\end{array}
\right)
\right\}
 \end{equation*}
and
$$
S=\mbox{Span}\left\{
\left(\begin{array}{l}
z_1^2z_2\\
0
\end{array}
\right),
\left(\begin{array}{l}
0\\
z_1z_2^2
\end{array}
\right)
\right\}.
$$
For notational convenience, in what follows we let
$$\mathcal{H}\left(\alpha z_1^{q_1}z_2^{q_2} \mu \right)=\left(\begin{array}{c}
\alpha z_1^{q_1}z_2^{q_2}\mu \vspace{0.2cm}\\
\overline{\alpha}z_1^{q_2}z_2^{q_1}\mu
\end{array}\right), ~\alpha\in \mathbb{C}.$$
We then calculate $\bf{g_j^1(z,0,\mu)}$.

\vspace{0.2cm}
\begin{subsubsection}
{Calculation of $\bf{g_2^1(z,0,\mu)}$}
\end{subsubsection}

From \eqref{FD}, we have
\begin{equation}
\label{FD2MU}
F^d_2 (\varphi, \mu)=F^d_{20} (\varphi)+\mu F^d_{21} (\varphi) ,
\end{equation}
and
\begin{equation}
\label{FD3MU}
F^d_3 (\varphi, \mu)=\mu F^d_{31}(\varphi), ~F^d_j (\varphi, \mu)=(0, 0)^T,  ~j=4, 5, \cdots.
\end{equation}
where
\begin{equation}
\label{FD23}
\left\{
\begin{array}{l}
F^d_{20} (\varphi)=-2d_{21}\tau_c\left(\begin{array}{c}
0\vspace{0.2cm}\\
\varphi^{(1)}_{x}(-1)\varphi^{(2)}_{x}(0)+\varphi^{(1)}_{xx}(-1)\varphi^{(2)}(0) \end{array}\right),\vspace{0.2cm}\\
F^d_{21} (\varphi)=2D_1\varphi_{xx}(0)+2D_2\varphi_{xx}(-1),\vspace{0.2cm}\\
F^d_{31} (\varphi)=-6d_{21}\left(\begin{array}{c}
0\vspace{0.2cm}\\
\varphi^{(1)}_{x}(-1)\varphi^{(2)}_{x}(0)+\varphi^{(1)}_{xx}(-1)\varphi^{(2)}(0)\end{array}\right),   \vspace{0.2cm}\\
\end{array}
\right.
\end{equation}

%
%

It is easy to verify that
\begin{equation}
\label{G211}
\left(\begin{array}{c}\left[2\mu A\left(\Phi(0) z_x\right),\beta_{n_c}^{(1)}\right]\vspace{0.3cm}\\
\left[2\mu A\left(\Phi(0) z_x\right),\beta_{n_c}^{(2)}\right]
\end{array}\right)=
2\mu  A\left( \Phi (0) \left(\begin{array}{c}z_1\\z_2\end{array}\right) \right),
\end{equation}
\begin{equation}
\label{G212}
\begin{array}{lll}
&&\left(\begin{array}{c}\left[\mu F^d_{21} \left(\Phi(\theta) z_x\right),\beta_{n_c}^{(1)}\right]\vspace{0.3cm}\\
\left[\mu F^d_{21} \left(\Phi(\theta) z_x\right),\beta_{n_c}^{(2)}\right]
\end{array}\right)\vspace{0.3cm}\\
&=&
-2(n_c/\ell)^2\mu \left(D_1\left( \Phi (0) \left(\begin{array}{c}z_1\\z_2\end{array}\right) \right)+D_2\left( \Phi (-1) \left(\begin{array}{c}z_1\\z_2\end{array}\right) \right)\right).
\end{array}
\end{equation}
In addition, from \eqref{FPM}, we have
for all $\mu\in \mathbb{R}$,
\begin{equation*}
\label{F2PHI}
F_2( \Phi(\theta) z_x,\mu)=F_2( \Phi(\theta) z_x, 0),
\end{equation*}

 It follows from  \eqref{FJ1} that
 \begin{equation}
\label{G2F21}
 f^1_2(z, 0,  \mu)=\Psi(0)\left(\begin{array}{c}\left[\widetilde{F}_2\left(\Phi(\theta) z_x,\mu\right),\beta_{n_c}^{(1)}\right]\vspace{0.1cm}\\
\left[\widetilde{F}_2\left(\Phi(\theta) z_x,\mu\right),\beta_{n_c}^{(2)}\right]
\end{array}\right).
\end{equation}
   This,  together with  \eqref{TLDF2MU},  \eqref{FD2MU}, \eqref{G211}  and  \eqref{G212}, yields to
\begin{equation}
\label{G21}
g^1_2(z,0,\mu)=\mbox{Proj}_{\tiny{\mbox{Ker}(M_2^1)}}f_2^1(z,0,\mu)
=
\mathcal{H}\left(B_1 \mu z_1\right),
\end{equation}
 where
$$
B_1=2\psi^{T}\left(A \phi (0) -(n_c/\ell)^2   \left(D_1\phi (0)+D_2\phi(-1)\right)\right)=2i\omega_{n_c}\psi ^{T}\phi. $$
\begin{subsubsection}
{Calculation of $\bf{g_3^1(z,0,\mu)}$}
\end{subsubsection}

In this subsection, we  calculate the
third term $g^1_3(z,0,0)$ in terms of \eqref{G310}. Notice that
  $\frac{1}{3!}\widetilde{f}_3^1$ in \eqref{G310} is the term of order 3
obtained after the changes of variables in previous step.
Denote
\begin{equation}
\label{F211W}
 f^{(1,1)}_2(z, w, 0)=\Psi(0)\left(\begin{array}{c}\left[F_2\left(\Phi(\theta) z_x+w, 0\right),\beta_{n_c}^{(1)}\right]\vspace{0.2cm}\\
\left[F_2\left(\Phi(\theta) z_x+w, 0\right),\beta_{n_c}^{(2)}\right]
\end{array}\right),
\end{equation}
\begin{equation}
\label{F212W}
 f^{(1,2)}_2(z, w, 0)=\Psi(0)\left(\begin{array}{c}\left[F^d_2\left(\Phi(\theta) z_x+w, 0\right),\beta_{n_c}^{(1)}\right]\vspace{0.2cm}\\
\left[F^d_2\left(\Phi(\theta) z_x+w, 0\right),\beta_{n_c}^{(2)}\right]
\end{array}\right).
\end{equation}

In addition, it follows from \eqref{G21}   that $g_2^1(z,0,0)=(0, 0)^T$. Then  $\widetilde{f}_3^1(z,0,0)$ is determined  by
\begin{equation}
\label{FTLD31}
\begin{array}{lll}
&&\widetilde{f}_3^1(z,0,0)\vspace{0.3cm}\\
&=&f_3^1(z,0,0)+\frac{3}{2}\left[\left(D_zf_2^1(z,0,0)\right)U_2^1(z,0)+\left(D_wf_2^{(1,1)}(z,0,0)\right)U_2^2(z,0)(\theta)\right.\vspace{0.3cm}\\
&&\left.+\left(D_{w,w_x,w_{xx}}f^{(1,2)}_2(z,0,0) \right)U_2^{(2,d)}(z, 0)(\theta)\right],
\end{array}
 \end{equation}
 where $f_2^1(z,0,0)=f_2^{(1,1)}(z,0,0)+f_2^{(1,2)}(z,0,0)$,
 $$D_{w,w_x,w_{xx}}f^{(1,2)}_2(z,0,0)=\left(D_wf^{(1,2)}_2(z,0,0), D_{w_x}f^{(1,2)}_2(z,0,0), D_{w_{xx}}f^{(1,2)}_2(z,0,0) \right),$$
  \begin{equation}
  \label{U212}
U_2^1(z,0)=\left(M_2^1\right)^{-1}\mbox{Proj}_{Im\left(M_2^1\right)}f_2^1(z,0,0),~~U_2^2(z,0)(\theta)=\left(M_2^2\right)^{-1}f_2^2(z,0,0),
 \end{equation}
 and
\begin{equation}
\label{U2D}
 U_2^{(2,d)}(z, 0)(\theta)=\mbox{col}\left(U_2^2(z, 0)(\theta), U_{2x}^2(z, 0)(\theta), U_{2xx}^2(z, 0)(\theta)\right).
 \end{equation}

 Next, we compute  $\mbox{Proj}_{S}\widetilde{f}_3^1(z,0,0)$ step by step according to \eqref{FTLD31}. The calculation is divided into the following four steps.

\vspace{0.4cm}
\noindent{\bf
 {Step 1: The calculation of $\bf{\mbox{Proj}_{S}f_3^1(z,0,0)}$ }}

Let
\begin{equation}
\label{F3EXP}
F_3(\Phi(\theta)
z_x,0)=\gamma_{n_c}^3(x) \left(\sum_{q_1+q_2=3}A_{q_1q_2}z_1^{q_1}z_2^{q_2}\right), ~q_1,q_2\in \mathbb{N}_0.
\end{equation}
 
From  \eqref{TLDF3MU} and \eqref{FD3MU},  we have
$\widetilde{F}_3(\Phi(\theta)z_x,0)= F_3(\Phi(\theta)
z_x,0)$.  Then it follows from \eqref{FJ1} and \eqref{F3EXP}  that
$$
f_3^1(z, 0, 0)=\Psi (0)\left(
\sum\limits_{q_1+q_2=3}A_{q_1q_2}\int_0^{\ell\pi}\gamma_{n_c}^{q_1+q_2+1}(x)dxz_1^{q_1}z_2^{q_2}
\right),
$$
which, together with
the fact that
$ \int_0^{\ell\pi}\gamma_{n_c}^4(x)  dx=
\frac{3}{2\ell\pi},
  $
implies that
\begin{equation}\label{PF31}
\mbox{Proj}_{S}f_3^1(z,0,0)=
\mathcal{H} \left(B_{21}z_1^2z_2\right)
\end{equation}
where
\begin{equation}
\label{B21}
B_{21}=\frac{3}{2\ell\pi}  \psi ^{T}A_{21} .
\end{equation}
\vspace{0.4cm}
\noindent{\bf Step 2: The calculation of $\bf{\mbox{Proj}_{S}\left(\left(D_zf_2^1\right)(z,0,0)U_2^1(z,0)\right)}$}

Form  \eqref{TLDF2MU} and   \eqref{FD2MU}, we have
\begin{equation}
\label{TLDF2}
\widetilde{F}_2( \Phi(\theta) z_x,0)=F_2( \Phi(\theta) z_x,0)+F^d_{20}\left (\Phi(\theta) z_x\right)
\end{equation}
By \eqref{FPM},  we  write
\begin{equation}
\label{F2EXP}
\begin{array}{lll}
&&F_2(\Phi(\theta) z_x+w,\mu) =F_2(\Phi(\theta)
z_x+w,0)   \vspace{0.2cm}\\
&=&\gamma_{n_c}^2(x)\left(\sum\limits_{q_1+q_2=2}A_{q_1q_2}z_1^{q_1}z_2^{q_2}\right)
+\mathcal{S}_2(\Phi(\theta)
z_x,w)+O\left(|w|^2\right), ~q_1,q_2 \in \mathbb{N}_0,
\end{array}
\end{equation}
where   $\mathcal{S}_2(\Phi(\theta)
z_x,w)$  is the second cross terms of $\Phi (\theta) z_x$ and $w$.
In addition,  by  \eqref{FD2MU} and \eqref{FD23},   we write
\begin{equation}
\label{F2DEXP}
\begin{array}{lll}
&& F^d_2\left (\Phi(\theta) z_x, 0\right)=F^d_{20}\left (\Phi(\theta) z_x\right) \vspace{0.2cm}\\
&=&(n_c/\ell)^2\left(\xi^2_{n_c}(x)- \gamma_{n_c}^2(x)\right) \left(\sum\limits_{q_1+q_2=2}A^d_{q_1q_2} z_1^{q_1}z_2^{q_2}\right),
\end{array}
\end{equation}
where
$$\xi_{n_c}(x)=\frac{\sqrt{2}}{\sqrt{\ell \pi}}\sin\left(\frac{n_cx}{\ell}\right),$$
\begin{equation}
\label{AIJD}
\left\{
\begin{array}{l}
 A_{20}^d=-
2d_{21}\tau_c\left(\begin{array}{c}
0\vspace{0.2cm}\\
\phi_1(-1) \phi_2(0)
\end{array}\right) =\overline{A_{02}^d},
 \vspace{0.3cm}\\
 A_{11}^d=-
2d_{21}\tau_c\left(\begin{array}{c}
0\vspace{0.2cm}\\
2\mbox{Re}\left\{\phi _1(-1) \overline{\phi_2}(0)\right\}
\end{array}\right).
\end{array}
\right.
\end{equation}

It is easy to verify that
$$
\int_0^{\ell\pi}\gamma_{n_c}^3(x)dx=\int_0^{\ell\pi}\xi_{n_2}^2(x)\gamma_{n_c}(x)dx =0,
$$
  Then, by \eqref{TLDF2}-\eqref{F2DEXP}, we have
$$
 f^1_2(z, 0, 0) = \Psi(0)\left(\begin{array}{c}\left[\widetilde{F}_2\left(\Phi(\theta) z_x, 0\right),\beta_{n_c}^{(1)}\right] \vspace{0.2cm}\\
 \left[\widetilde{F}_2\left(\Phi(\theta) z_x, 0\right),\beta_{n_c}^{(2)}\right]
 \end{array}
 \right)
 = \left(\begin{array}{c}\
 0\\
 0
 \end{array}\right)
$$

Hence,
\begin{equation}
\label{PDZF21}
\mbox{Proj}_{S}\left[\left(D_zf_2^1\right)(z,0,0)U_2^1(z,0)\right]= \left(\begin{array}{c}\
 0\\
 0
 \end{array}\right).
\end{equation}
\vspace{0.5cm}
\noindent{\bf
 Step 3: The calculation of $\bf{\mbox{Proj}_{S}\left(\left(D_wf_2^{(1,1)}(z,0,0)\right)U_2^2(z,0)(\theta)\right)}$}

 Let
\begin{equation}
\label{U22}
 U_2^2(z, 0)(\theta)\triangleq h(\theta, z)=\sum\limits_{n\in \mathbb{N}_0}  h_n(\theta,z) \gamma_n(x),
 \end{equation}
where
$$  h_n(\theta,z) =\sum\limits_{q_1+q_2=2} h_{n,q_1q_2}(\theta)z_1^{q_1}z_2^{q_2}.$$

Then, from \eqref{F211W} and \eqref{U22},  we obtain
 $$
\begin{array}{lll}
&& \left(D_wf^{(1,1)}_2(z,0,0)\right)U_2^2(z,0)(\theta)\vspace{0.2cm}\\
&=&\Psi(0)
   \left(\begin{array}{c}
\left[\left.D_w F_2\left(\Phi(\theta) z_x+w, 0\right)\right|_{w=0}\left(\sum\limits_{n\in \mathbb{N}_0}  h_n(\theta,z) \gamma_n(x)\right), \beta_{n_c}^{(1)}\right]\vspace{0.2cm}\\
\left[\left.D_w F_2\left(\Phi(\theta) z_x+w, 0\right)\right|_{w=0}\left(\sum\limits_{n\in \mathbb{N}_0}  h_n(\theta,z) \gamma_n(x)\right), \beta_{n_c}^{(2)}\right]
\end{array}\right).
\end{array}
 $$

By \eqref{F2EXP}, we obtain
$$\left.D_w F_2\left(\Phi(\theta) z_x+w, 0\right)\right|_{w=0}\left(\sum\limits_{n\in \mathbb{N}_0}  h_n(\theta,z) \gamma_n(x)\right)= \mathcal{S}_2\left(\Phi(\theta) z_x, \sum\limits_{n\in \mathbb{N}_0}  h_n(\theta,z) \gamma_n(x)\right) $$
and
$$
\begin{array}{lll}
&&\left(\begin{array}{c}
                        \left[\mathcal{S}_2\left(\Phi(\theta) z_x, \sum\limits_{n\in \mathbb{N}_0}  h_n(\theta,z) \gamma_n(x)\right), \beta_{n_c}^{(1)}\right] \vspace{0.2cm}\\
  \left[\mathcal{S}_2\left(\Phi(\theta)  z_x, \sum\limits_{n\in \mathbb{N}_0}  h_n(\theta,z) \gamma_n(x)\right), \beta_{n_c}^{(2)}\right]
                         \end{array}
\right)\vspace{0.2cm}\\
&=&\sum\limits_{n\in \mathbb{N}_0} b_n\left(\mathcal{S}_2\left(\phi (\theta) z_1,h_n(\theta, z)\right)+\mathcal{S}_2\left(\overline{\phi }(\theta) z_2,h_n(\theta, z)\right)\right),
\end{array}
$$
where
\begin{equation}
\label{BN}
b_n=\int_0^{\ell\pi}\gamma^2_{n_c}(x)\gamma_n(x) dx=\left\{\begin{array}{ll}
\frac{1}{\sqrt{\ell\pi}}, & n=0,      \vspace{0.2cm}\\
\frac{1}{\sqrt{2\ell\pi}}, & n=2n_c,    \vspace{0.2cm}\\
 0, & \mbox{otherwise}.
 \end{array}
\right.
\end{equation}

Hence,
$$
\begin{array}{lll}
&&\left(D_wf_2^{(1,1)}(z, 0,
0)\right)U^2_2(z,0)(\theta)\vspace{0.2cm}\\
&=&\Psi(0)\left(
\begin{array}{l}
\sum\limits_{n=0, 2n_c} b_n\left(\mathcal{S}_2\left(\phi (\theta) z_1,h_n(\theta, z)\right)+\mathcal{S}_2\left(\overline{\phi }(\theta) z_2,h_n(\theta, z)\right)\right)
\end{array}
\right).
\end{array}
$$

Then, we have
\begin{equation}
\label{PDWF221}
\textmd{Proj}_{S}\left(D_wf_2^{(1,1)}(z, 0,
0)U^2_2(z,0)(\theta)\right)=
\mathcal{H} \left(B_{22}z_1^2z_2\right),
 \end{equation}
where
\begin{equation}
\label{B22}
\begin{array}{lll}
B_{22} &=&\frac{1}{\sqrt{\ell\pi}}\psi ^T\left(\mathcal{S}_2(\phi(\theta) ,h_{0,11}(\theta))+\mathcal{S}_2(\overline{\phi} (\theta) ,h_{0,20}(\theta) ) \right) \vspace{0.2cm}\\
&&+\frac{1}{\sqrt{2\ell\pi}}\psi ^T\left(\mathcal{S}_2(\phi(\theta) ,h_{2n_c, 11}(\theta))+\mathcal{S}_2(\overline{\phi} (\theta) ,h_{2n_c, 20}(\theta) ) \right).
\end{array}
\end{equation}
\vspace{0.4cm}
\noindent
{\bf Step 4: The calculation of $\bf{ \mbox{Proj}_{S}\left(\left(D_{w,w_x,w_{xx}}f^{(1,2)}_2(z,0,0)\right)U_2^{(2,d)}(z, 0)(\theta)\right)}$}

Denote $\varphi(\theta)=\left(\varphi^{(1)},\varphi^{(2)}\right)^T=  \Phi(\theta) z_x$ and
$$
\begin{array}{lll}
&&F_2^d( \varphi(\theta),  w, w_x, w_{xx})=F_2^d( \varphi(\theta)+w,0)= F_{20}^d( \varphi(\theta)+w)  \vspace{0.3cm}\\
&=&
-2d_{21}\tau_c\left(\begin{array}{c}
0\vspace{0.2cm}\\
\left(\varphi^{(1)}_{xx}(-1)+w_{xx}^{(1)}(-1)\right)\left(\varphi^{(2)}(0)+w^{(2)}(0)\right)
\end{array}\right)
\vspace{0.3cm}\\
&&
-2d_{21}\tau_c\left(\begin{array}{c}
0\vspace{0.2cm}\\
\left(\varphi^{(1)}_x(-1)+w_{x}^{(1)}(-1)\right)\left(\varphi^{(2)}_{x}(0)+w_{x}^{(2)}(0)\right)
\end{array}\right),
\end{array}
$$

and for $\phi(\theta)=\left(\phi_1(\theta), \phi_2(\theta)\right )^T, y(\theta)=\left(y_1(\theta), y_2(\theta)\right )^T\in C\left([-1, 0], \mathbb{R}^2\right)$,
\begin{equation}
\label{SDJ}
\left\{
\begin{array}{l}
\mathcal{S}_2^{(d,1)}(\phi (\theta), y(\theta))=
-2d_{21}\tau_c\left(\begin{array}{c}
0\vspace{0.2cm}\\
\phi_1(-1)y_2(0)
\end{array}\right),
\vspace{0.3cm}\\
\mathcal{S}_2^{(d,2)}(\phi (\theta), y(\theta) )
=
-2d_{21}\tau_c\left(\begin{array}{c}
0\vspace{0.2cm}\\
\phi_1(-1)y_2(0)
\end{array}\right)-2d_{21}\tau_c\left(\begin{array}{c}
0\vspace{0.2cm}\\
\phi_2(0)y_1(-1)
\end{array}\right),

\vspace{0.3cm}\\
\mathcal{S}_2^{(d,3)}(\phi(\theta), y(\theta) )=
-2d_{21}\tau_c\left(\begin{array}{c}
0\vspace{0.2cm}\\
\phi_2(0)y_1(-1)
\end{array}\right).
\end{array}
\right.
\end{equation}

Then, from \eqref{F212W},  \eqref{U2D} and \eqref{U22},   we have
       $$
\begin{array}{lll}
&& \left(D_{w,w_x,w_{xx}}f^{(1,2)}_2(z,0,0)\right)U_2^{(2,d)}(z, 0)(\theta)\vspace{0.3cm}\\
 &=&\Psi(0)
   \left(\begin{array}{c}
\left[D_{w, w_x, w_{xx}}F_2^d( \varphi(\theta),  w, w_x, w_{xx}) U_2^{(2,d)}(z, 0)(\theta), \beta_{n_c}^{(1)}\right]\vspace{0.2cm}\\
\left[D_{w, w_x, w_{xx}}F_2^d( \varphi(\theta),  w, w_x, w_{xx})U_2^{(2,d)}(z, 0)(\theta), \beta_{n_c}^{(2)}\right]
\end{array}\right)
\end{array}
$$
and then  we obtain
   \begin{equation}
\label{PDWF22}
\textmd{Proj}_{S}\left(\left(D_{w,w_x,w_{xx}}f^{(1,2)}_2(z,0,0)\right)U_2^{(2,d)}(z, 0)(\theta)\right)=
\mathcal{H} \left(B_{23}z_1^2z_2\right),
 \end{equation}
where
\begin{equation}
\label{B23}
\begin{array}{lll}
B_{23}
&=& -\frac{1}{\sqrt{\ell\pi}} \left( n_c/\ell\right)^2 \psi ^T\left(\mathcal{S}_2^{(d,1)}\left(\phi (\theta), h_{0,11}(\theta)\right)+\mathcal{S}_2^{(d,1)}\left(\overline{\phi} (\theta), h_{0,20}(\theta)\right)\right)
 \vspace{0.2cm}\\
&&+\frac{1}{\sqrt{2\ell\pi}}\psi ^T\sum\limits_{j=1,2,3} b_{2n_c}^{(j)}  \left( \mathcal{S}_2^{(d,j)}\left(\phi (\theta), h_{2n_c, 11}(\theta)\right) + \mathcal{S}_2^{(d,j)}\left(\overline{\phi}(\theta), h_{2n_c, 20}(\theta)\right)\right)
\end{array}
\end{equation}
with
$$b_{2n_c}^{(1)}=-(n_c/\ell)^2, ~b_{2n_c}^{(2)}=2(n_c/\ell)^2, ~ b_{2n_c}^{(3)}=-(2n_c/\ell)^2.$$

\vspace{0.2cm}
\begin{subsubsection}
{Normal form of the Hopf bifurcation and the corresponding coefficients}
\end{subsubsection}
 According to the above calculations, we have obtained the normal form of the Hopf bifurcation in the following form
\begin{equation}
\label{NF-Hopf}
  \dot{z}=Bz+\frac{1}{2}\left(
               \begin{array}{c}
                 B_1 z_1\mu \\
                 \bar{B}_1z_2\mu \\
               \end{array}
             \right)+\frac{1}{3!}\left(
                       \begin{array}{c}
                         B_2 z_1^2 z_2 \\
                         \bar{B}_2 z_1 z^2_2 \\
                       \end{array}
                     \right)+O(|z|\mu^2+|z^4|),
\end{equation}
where
$$B_1=2i\omega_{n_c}\psi^{T}\phi , ~B_2=B_{21}+\frac{3}{2}\left(B_{22}+B_{23} \right),$$
and $B_{2j}, j=1,2,3,$ are determined by  \eqref{B21}, \eqref{B22} and \eqref{B23}.
Through the change of
variables  $z_1=w_1-iw_2, z_2=w_1+iw_2$ and  $ w_1=\rho\cos\xi$, $w_2=\rho\sin\xi$, the normal form \eqref{NF-Hopf}  can be written as  the following form in polar
coordinates
$$
  \begin{array}{rcl}
   \dot{\rho} &=& K_1\mu\rho+K_2\rho^3+O(\mu^2\rho+|(\rho,\mu)|^4),\\
  \end{array}
  $$
with $K_1=\frac{1}{2}\textrm{Re}\left(B_1\right)$, $K_2=\frac{1}{3!}\textrm{Re}\left(B_2\right)$.
According to \cite{Chow-82}, the sign of $K_1K_2$ determines the direction
of the bifurcation (supercritical for $K_1K_2<0$  and subcritical for $K_1K_2>0$), and
  the sign of $K_2$  determines the stability of the Hopf bifurcating periodic orbits  (stable
for  $K_2<0$  and  unstable for $K_2>0$).


From \eqref{B21}, \eqref{B22} and \eqref{B23}, it is obvious that in order to obtain the value of $K_2$,
we still need to compute  $h_{0,20}(\theta), h_{0,11}(\theta), h_{2n_c, 20}(\theta), h_{2n_c, 11}(\theta)$, and $A_{ij}$.


From  \cite{Faria-00TAMS} , we have
$$
\begin{array}{lll}
&&M_2^2\left( h_n(\theta,z) \gamma_n(x)\right)\vspace{0.2cm}\\
&=& D_z\left( h_n(\theta,z) \gamma_n(x)\right)Bz -A_{\mathcal{Q}^1}\left( h_n(\theta,z) \gamma_n(x)\right),
\end{array}
$$
which leads to
 \begin{equation}
 \label{M22HJ}
 \begin{array}{lll}
 &&   \left(\begin{array}{c}
\left[M_2^2\left( h_n(\theta,z) \gamma_n(x)\right), \beta_n^{(1)}\right]\vspace{0.2cm}\\
\left[M_2^2\left( h_n(\theta,z) \gamma_n(x)\right), \beta_n^{(2)}\right]
\end{array}\right)\vspace{0.2cm}\\
&=&2i\omega_c\left(h_{n,20} (\theta)z_1^2- h_{n,02}(\theta)z_2^2\right)\vspace{0.2cm}\\
&&-\left(\dot{h}_n(\theta, z)+X_0(\theta)\left(\mathscr{L}_0 \left(h_n(\theta, z)\right)- \dot{h}_n(0, z)\right)\right),
\end{array}
\end{equation}
where
$$
\mathscr{L}_0\left( h_n(\theta, z)\right)= -\tau_c(n/\ell)^2\left( D_1 h_n(0, z)+D_2h_n(-1, z) \right)+\tau_cAh_n(0, z).
$$

By \eqref{FJ2}, we get
\begin{equation}
\label{F22}
\begin{array}{lll}
&&f_2^2(z,0,0)\vspace{0.2cm}\\
&=&X_0(\theta)\widetilde{F}_2\left(\Phi(\theta) z_x,0\right) -\pi\left(X_0(\theta) \widetilde{F}_2\left(\Phi(\theta) z_x,0\right)\right) \vspace{0.2cm}\\
&=&X_0(\theta)\widetilde{F}_2\left(\Phi(\theta) z_x,0\right) -\Phi (\theta)\Psi (0)\left(\begin{array}{c}\left[\widetilde{F}_2\left(\Phi(\theta) z_x,0\right),\beta_{n_c}^{(1)}\right]\vspace{0.2cm}\\
\left[\widetilde{F}_2\left(\Phi(\theta) z_x,0\right),\beta_{n_c}^{(2)}\right]
\end{array}\right)\gamma_{n_c}(x). \end{array}
\end{equation}
By \eqref{PIX0},  \eqref{TLDF2}, \eqref{F2EXP} and \eqref{F2DEXP}, we have
\begin{equation}\label{EQF22B}
\begin{array}{lll}
 &&
\left(\begin{array}{c}
\left[f_2^2(z,0,0),\beta_n^{(1)}\right]\vspace{0.1cm}\\
\left[f_2^2(z,0,0),\beta_n^{(2)}\right]
\end{array}\right)\vspace{0.2cm}\\
&=&\begin{cases}
\begin{array}{l}
\frac{1}{\sqrt{\ell\pi}} X_0(\theta)\left(A_{20}z_1^2+A_{02}z_2^2+A_{11}z_1z_2\right),
\end{array}
& n=0,\vspace{0.3cm}\\
\frac{1} {\sqrt{2\ell\pi}} X_0(\theta)\left(\widetilde{A}_{20}z_1^2+\widetilde{A}_{02}z_2^2+\widetilde{A}_{11}z_1z_2\right),
&n=2n_c,\vspace{0.3cm}\\
\end{cases}
\end{array}
\end{equation}
where $\widetilde{A}_{j_1j_2} $ is defined as follows
\begin{equation*}
\label{WLDAJ3}
\left\{\begin{array}{l}
\widetilde{A}_{j_1j_2}=A_{j_1j_2}- 2(n_c/\ell)^2A_{j_1j_2}^d,\vspace{0.2cm}\\
j_1, j_2=0, 1, 2, \quad j_1+j_2=2,
\end{array}\right.
\end{equation*}
where  $A_{j_1j_2}^d$ is determined by \eqref{AIJD} and $A_{j_1j_2}$ will be determined in the following.

Hence, from \eqref{U212}, \eqref{M22HJ}-\eqref{EQF22B} and matching the
coefficients of $z_1^2, z_1z_2 $, we have
\begin{equation}
\label{H0_EQ}
n=0,~\left\{\begin{array}{ll}
z_1^2: &\begin{cases}
\dot{h}_{0,20}(\theta)-2i\omega_c h_{0,20}(\theta)=(0, 0)^T,\vspace{0.1cm}\\
\dot{h}_{0,20}(0)-L_0(h_{0,20}(\theta))=\frac{1}{\sqrt{\ell\pi}}A_{20},
\end{cases}
\vspace{0.2cm}\\
z_1z_2:&\begin{cases}
\dot{h}_{0,11}(\theta)=(0, 0)^T,\vspace{0.1cm}\\
\dot{h}_{0,11}(0)-L_0(h_{0,11}(\theta))=\frac{1}{\sqrt{\ell\pi}}A_{11},
\end{cases}
\end{array}
\right.
\end{equation}
\begin{equation}
\label{H2N11_EQ}
n=2n_c,
~\left\{\begin{array}{ll}
z_1^2: &\begin{cases}
\dot{h}_{2n_c,20}(\theta)-2i\omega_c h_{2n_c,20}(\theta)=(0, 0)^T,\vspace{0.1cm}\\
\dot{h}_{2n_c,20}(0)-\mathscr{L}_0(h_{2n_c,20}(\theta))=\frac{1}{\sqrt{2\ell\pi}}\widetilde{A}_{20},
\end{cases}\vspace{0.2cm}\\
z_1z_2:&\begin{cases}
\dot{h}_{2n_c,11}(\theta)=(0, 0)^T,\vspace{0.1cm}\\
\dot{h}_{2n_c,11}(0)-\mathscr{L}_0(h_{2n_c,11}(\theta))=\frac{1}{\sqrt{2\ell\pi}}\widetilde{A}_{11},
\end{cases}
\end{array}
\right.
\end{equation}
Solving \eqref{H0_EQ}, \eqref{H2N11_EQ}, we obtain
\begin{equation*}
\label{H0}
\left\{
\begin{array}{l}
h_{0,20}(\theta) =\frac{1}{\sqrt{\ell\pi}}\left(\widetilde{\mathcal{M}}_0(2i\omega_c)\right)^{-1}A_{20}e^{2i\omega_c\theta},\vspace{0.2cm}\\
h_{0,11}(\theta) =\frac{1}{\sqrt{\ell\pi}}\left(\widetilde{\mathcal{M}}_0(0)\right)^{-1}A_{11},
\end{array}
\right.
\end{equation*}
and
\begin{equation*}
\label{H2NC}
\left\{
\begin{array}{lll}
h_{2n_c,20}(\theta) &=&\frac{1}{\sqrt{2\ell\pi}}\left(\widetilde{\mathcal{M}}_{2n_c}(2i\omega_c)\right)^{-1}\widetilde{A}_{20}e^{2i\omega_c\theta}\vspace{0.2cm}\\
h_{2n_c,11}(\theta) &=&\frac{1}{\sqrt{2\ell\pi}}\left(\widetilde{\mathcal{M}}_{2n_c}(0)\right)^{-1}\widetilde{A}_{11},
\end{array}
\right.
\end{equation*}
where the matrix $ \widetilde{\mathcal{M}}_n(\lambda)$ is defined by  \eqref{NFMN}.

\begin{remark}{\it
Although the normal form of the delay-induced Hopf bifurcation has extensively studied, the existing procedures can not be applied to the memory-based diffusion system because of the
nonlinearity of the diffusion terms and the existence of the delay in the diffusion terms.  Compared with the existing procedures, the procedure in our Step 4 is the main characteristic difference.}
\end{remark}
\section{Application to a  predator-prey model with Holling type-II functional response}\label{sec3}

 We now apply our newly developed algorithm in Section \ref{sec2} to a a  predator-prey model, where in the model \eqref{PPMODEL}
 $$
f(u,v)=u\left(1-\frac{u}{a}\right)-\frac{buv}{1+u},~~g(u,v)=\frac{buv}{1+u}-c v.
$$
 Then,  \eqref{PPMODEL}  becomes the following  predator-prey model with Holling type II functional response:
\begin{equation}
\label{Exa-PP}
\left\{  \begin{array}{ll}
  \frac{\partial u(x,t)}{\partial t}=d_{11}u_{xx}(x,t)+u\left(1-\frac{u}{a}\right)-\frac{buv}{1+u}, & 0<x<\ell\pi, t>0, \vspace{0.2cm}\\
\frac{\partial v(x,t)}{\partial t}=-d_{21}(v(x,t)u_{x}(x,t-\tau))_x+d_{22}v_{xx}(x,t)-cv+\frac{buv}{1+u}, & 0<x<\ell\pi, t>0,\vspace{0.2cm}\\
u_x(0,t)=u_x(\ell\pi,t)=v_x(0,t)=v_x(\ell\pi,t)=0, & t\geq 0.
  \end{array}
  \right.
 \end{equation}

\begin{subsection}
   {Stability and bifurcation analysis}
 \end{subsection}

System \eqref{Exa-PP} has the  positive constant steady state $E_*(u_*, v_*)$, where
$$u_*=\frac{c}{b-c},~~v_*=\frac{(a-u_*)(1+u_*)}{ab},$$
provided that $b>\frac{c(1+a)}{a}$ holds.  For  $E_*(u_*, v_*)$,  we have
\begin{equation}
\label{AIJ-EXA}
\begin{array}{c}
a_{11}=\frac{\gamma(a-1-2\gamma)}{a(1+\gamma)}
\left\{\begin{array}{ll}
\leq 0, & \frac{a-1}{2}\leq \gamma< a,\vspace{0.2cm}\\
>0, & 0<\gamma<\frac{a-1}{2},
\end{array}\right.
 \vspace{0.2cm}\\
 a_{12}=-c<0, ~~a_{21}=\frac{a-\gamma}{a(1+\gamma)}>0,   ~~a_{22}=0,
\end{array}
\end{equation}
where $\gamma=\frac{c}{b-c}$.
Let
\begin{equation}
\label{JN}
J_n=d_{11}d_{22} (n/\ell)^4
-\left(d_{11}a_{22}+d_{22}a_{11}\right) (n/\ell)^2+Det(A),
\end{equation}
From \eqref{AIJ-EXA} and \eqref{JN}, it is easy to verify that $T_n<0$ and $J_n>0$  provided that
$$(C_0)\quad \quad \frac{a-1}{2}< \gamma< a.$$
This implies that   when $d_{21}=0$ and the condition $(C_0)$ holds,  the  positive steady state $E_*$ is asymptotically  stable for   $d_{11}\geq 0$ and $d_{22}\geq 0$
 In what follows, we always assume that the condition $(C_0)$ holds.

  Since $J_n>0 $ under the condition $(C_0)$, we have  $\Gamma(0)=J_n-d_{21}v_*a_{12}>0$.    This implies that $\lambda=0$ is not a root of Eq.\eqref{CE}.

Let  $\lambda=i \omega~(\omega>0)$ be a root of \eqref{CE}.  Substituting it  along with expressions in \eqref{AIJ-EXA} into \eqref{CE} and separating the real from the imaginary parts, we have

\begin{equation}
\label{RIPARTS}
  \begin{cases}
\omega^2-J_n= (n/\ell)^2cd_{21}v_*\cos\left(\omega\tau\right),\vspace{0.1cm}\\
-T_n\omega=(n/\ell)^2cd_{21}v_*\sin\left(\omega\tau\right),
  \end{cases}
 \end{equation}
which yields
\begin{equation}
\label{OMGE}
\omega^4+P_n\omega^2+Q_n=0,
\end{equation}
where
\begin{equation}
\label{PN}
\begin{array}{lll}
P_n&=&T_n^2-2J_n\vspace{0.2cm}\\
&=&\left(d_{11}^2+d_{22}^2\right)(n/\ell)^4-2\left(d_{11}a_{11}+d_{22}a_{22}\right)(n/\ell)^2
+a_{11}^2+a_{22}^2+2a_{12}a_{21}.
\end{array}
\end{equation}
and
\begin{equation}
\label{QN}
Q_n=\left(J_n+cd_{21}v_*(n/\ell)^2\right)\left(J_n-cd_{21}v_*(n/\ell)^2\right).
\end{equation}

and
$$
a_{11}^2+a_{22}^2+2a_{12}a_{21}\left\{\begin{array}{ll}
\leq 0, &  c\geq c_*,\vspace{0.2cm}\\
>0, &  c<c_*,
\end{array}\right.
$$
where
\begin{equation}
\label{C-star}
c_*=\frac{\gamma^2(a-1-2\gamma)^2}{2a(1+\gamma)(a-\gamma)}.
\end{equation}

For simplification, we assume that $c<c_*$. Then, $P_n>0$ for any $n\in\mathbb{N}_0$.
Define
\begin{equation}
\label{d21n}
d_{21}^{(n)}  =\frac{J_n}{cv_*(n/\ell)^2}=\frac{1}{cv_*}\left( d_{11}d_{22} (n/\ell)^2 +\frac{Det(A)}{(n/\ell)^2}-\left( d_{11}a_{22}+d_{22}a_{11} \right)\right)>0,
\end{equation}
Then, for fixed $n$,  by \eqref{QN} we have
\begin{equation}
\label{QN1}
 Q_n\left\{
 \begin{array}{ll}
 >0, & 0<d_{21}<d_{21}^{(n)},   \vspace{0.2cm}\\
  =0, & d_{21}=d_{21}^{(n)},  \vspace{0.2cm}\\
  <0, & d_{21}>d_{21}^{(n)}.
 \end{array}
 \right.
\end{equation}
Thus, when $d_{21}>d_{21}^{(n)}$, Eq.\eqref{OMGE} has one positive root  $\omega_n$, where
\begin{equation}
\label{OMGN}
\omega_n =\sqrt{\frac{-P_n\pm \sqrt{P_n^2-4Q_n}}{2}}.
\end{equation}
Notice that $T_n<0$. From \eqref{RIPARTS}, we have
$$\sin\left(\omega_n\tau\right)=\frac{-T_n\omega_n}{(n/\ell)^2cd_{21}v_*}>0$$
Thus, set
\begin{equation}
\label{TAUNJ}
\tau_{n, j}=\frac{1}{\omega_n}\left\{\arccos \left\{\frac{\omega_n^2-J_n}{d_{21}v_*c(n/\ell)^2}\right\}+2j\pi\right\}, ~~j\in\mathbb{N}_0,~n\in\mathbb{N},
\end{equation}
 then \eqref{CE} has a pair of purely imaginary roots $\pm \omega_ni$ at $\tau=\tau_{n, j}^+$.  And it is easy to verify
the  transversality condition satisfies
 $\left.\frac{d \mathrm{Re}(\lambda (\tau))}{d
\tau}\right|_{\tau=\tau_{n,j}}>0$.

Let
 \begin{equation}
\label{D21STAR}
d_{21}^{*}=\min\limits_{n\in \mathbb{N}}\left\{ d_{21}^{(n)} \right\}>0.
\end{equation}
From \eqref{d21n}, it is easy to verify that $d_{21}^{~\left(n\right)}$ is decreasing for $n<\ell \sqrt[4]{\frac{Det(A)}{d_{11}d_{22}} }$, is increasing for $n>\ell\sqrt[4]{\frac{Det(A)}{d_{11}d_{22}} }$ and $d_{21}^{~\left(n\right)}\to\infty$ as $n\to\infty$.   This implies that
   $d_{21}^{*}$ exists.

 For fixed $d_{21}>d_{21}^*$, define an index set
  $$U(d_{21})=\left\{n\in \mathbb{N}: ~d_{21}^{(n)}  <d_{21}\right\}. $$
By the above discussion, we obtain the following results on the stability and Hopf bifurcation of system \eqref{Exa-PP}. 
 \begin{theorem}
 \label{The-31}
 Assume that  the condition $(C_0)$ holds and $c<c_*$.    Then we have the following:
    \begin{description}
    \item[(a)] when $0<d_{21}\leq d_{21}^*$,  the positive constant steady state $(u_*, v_*)$  of system \eqref{Exa-PP}  is locally  asymptotically stable   for any $\tau\geq 0$;
\item[(b)] when  $d_{21}>d_{21}^*$, there exists a critical value $\tau_*(d_{21})$ of the delay such that
\begin{itemize}
\item[(b1)] the positive constant steady state $(u_*, v_*)$  of system \eqref{Exa-PP}  is locally  asymptotically stable   for   $0\leq \tau<\tau_*(d_{21})$ and unstable for $ \tau>\tau_*(d_{21})$, where
$$\tau_*(d_{21})=\min\limits_{n\in U(d_{21})}\left\{ \tau_{n, 0}\right\}.$$
\item[(b2)]   system \eqref{Exa-PP}  undergoes mode$-n$ Hopf bifurcation at $\tau=\tau_{n,j}, n\in U(d_{21})$, and multiple Hopf bifurcations occur when these  Hopf bifurcation curves $\tau=\tau_{n,j}, n\in U(d_{21})$  have  interaction point on the $d_{21}-\tau$ plane.
\end{itemize}
\end{description}
  \end{theorem}

\begin{subsection}
{Direction and stability of the Hopf bifurcation}
\end{subsection}

We now numerically investigate the bifurcation direction and stability of the bifurcation. To this end, we set the parameters as follows
$$a=1, ~b=\frac{3}{10}, ~c=\frac{1}{10}, ~d_{11}=\frac{3}{5}, ~d_{22}=\frac{4}{5},~\ell=2.$$
Then, we have $(u_*, v_*)=(1/5, 5/2)$, $\gamma=1/2$,
$$
a_{11} =-\frac{1}{3}, a_{12} =-\frac{1}{10}, a_{21} =\frac{1}{3},  a_{22}=0
$$
and $Tr(A)=-1/2, Det(A)=1/30$.
  It follows from \eqref{PN}  and \eqref{d21n}  that
$$P_n=\frac{1}{16}n ^4 +\frac{1}{10}n^2+\frac{2}{45}>0,  $$
and
\begin{equation}
\label{NSd21N}
d_{21}^{(n)}=\frac{12n^2}{25} +\frac{8}{15n^2} + \frac{16}{15}.
\end{equation}

\begin{figure}
\centering
{\includegraphics[scale=0.5]{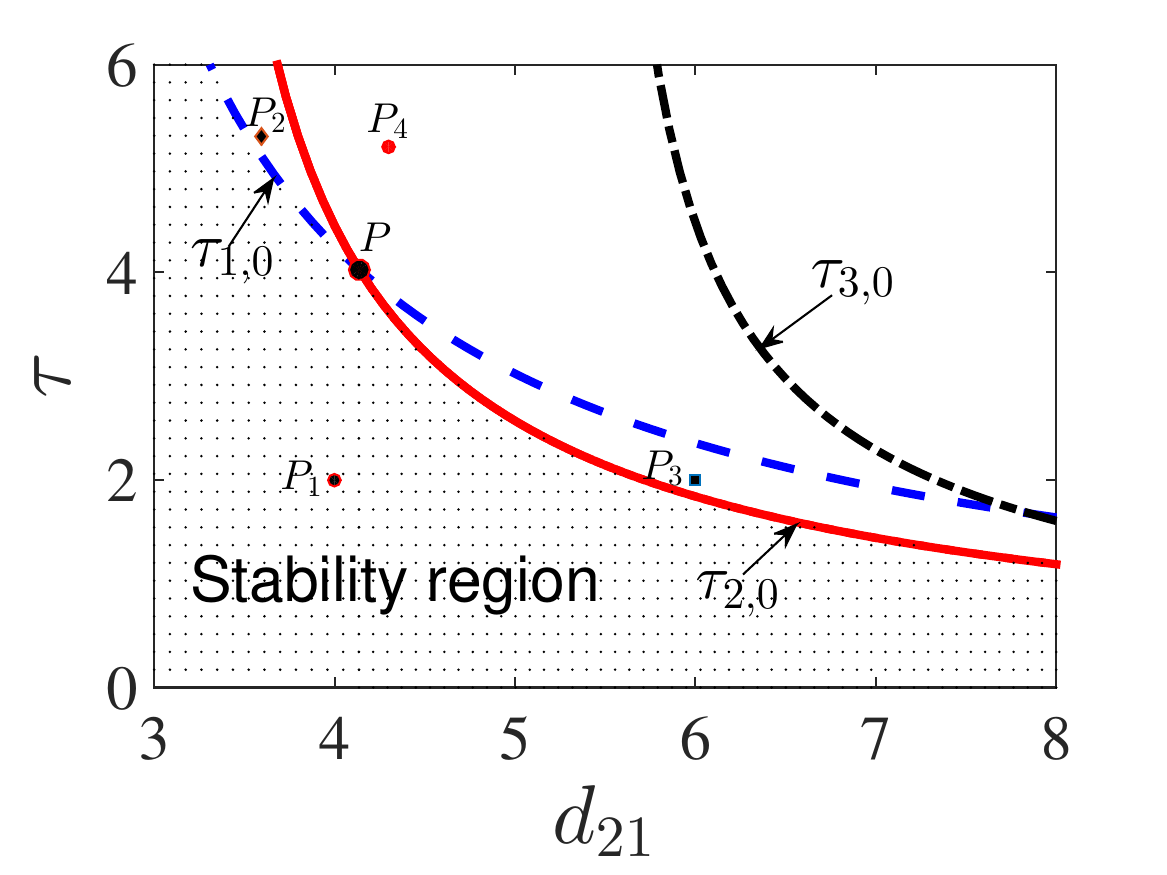}}
 \caption{Stability region and Hopf bifurcation curves in $d_{21}-\tau$ plane. The dotted region is the stability region and $\tau=\tau_{k, 0}, k=1,2,3$, are Hopf bifurcation curves. Hopf bifurcation curves $\tau=\tau_{1,0}$  and $\tau=\tau_{2,0}$ intersect at the point $P(4.1354,4.0292)$. The points $P_1(4, 2)$, $P_2(3.6,  5.3)$, $P_3(6,  2)$ and $P_4(4.3, 5.2)$ are chosen for the numerical simulations.  }
 \label{Fig1}
\end{figure}

\begin{figure}
\centering
\includegraphics[scale=0.3]{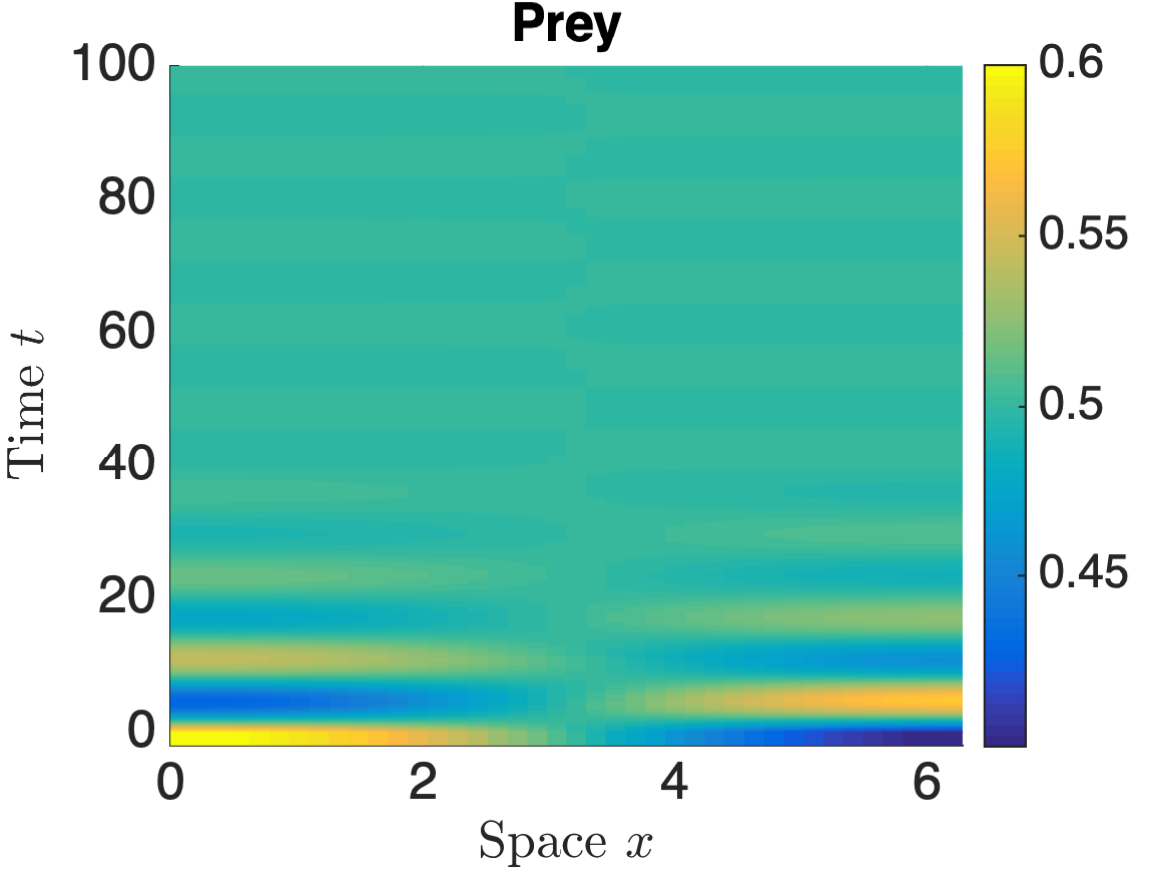}
\includegraphics[scale=0.3]{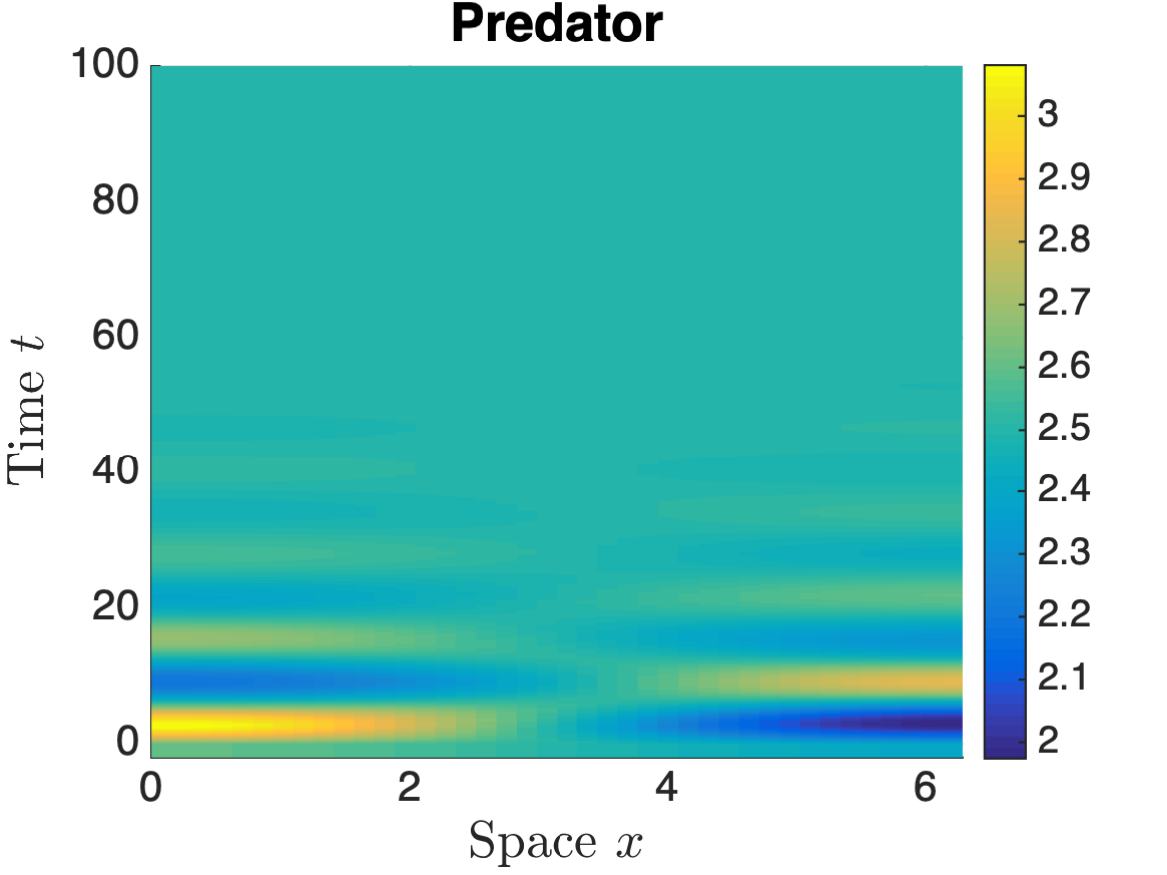}
  \caption{The spatial-temporal dynamics of system \eqref{Exa-PP}  with the parameters $a=1, ~b=\frac{3}{10}, ~c=\frac{1}{10}, ~d_{11}=\frac{6}{10}, ~d_{22}=\frac{8}{10},~\ell=2$ and   $\left(d_{21}, \tau\right)$ being chosen as the  point $P_1(4, 2)$ of Fig.\ref{Fig1}. The constant steady state is stable is asymptotically stable.  The initial values are $u_0(x)=0.2+0.1\cos(x/2),~v_0(x)=2.5+0.1\cos(x/2)$.}
  \label{Fig2}
  \end{figure}

Notice that $P_n>0$ for any $n\in \mathbb{N}$, which together with \eqref{QN1}, implies that for a fixed $n$,
 Eq.\eqref{OMGE} has no positive root for
$d_{21}<d_{21}^{(n)}$ and has only one positive root for $d_{21}\geq d_{21}^{(n)}$.
From \eqref{NSd21N},  it is easy to verify that $d_{21}^{(n)}<d_{21}^{(n+1)}$  for any $n\in \mathbb{N}$, and
$$
d_{21}^{(1)}=2.08< d_{21}^{(2)}=3.12<d_{21}^{(3)} \doteq 5.4459 .
$$
 Therefore,  by \eqref{D21STAR}, we have    $d_{21}^*=d_{21}^{(1)}=2.08$.   It follows from \eqref{C-star} that $c_*=1/6$. By Theorem  \eqref{The-31},  we  have the following stability result.
 \begin{proposition}
 \label{Pro-1}
For system \eqref{Exa-PP}  with the parameters $a=1, ~b=\frac{3}{10}, ~c=\frac{1}{10}, ~d_{11}=\frac{6}{10}, ~d_{22}=\frac{8}{10},~\ell=2$,
when $0\leq d_{21} < d_{21}^{(1)}=2.08$,   the positive constant steady state $(u_*, v_*)=(1/5,  5/2)$     is locally  asymptotically stable   for any $\tau\geq 0$;
  \end{proposition}
For fixed  $d_{21}>d_{21}^{(1)}$,    the positive constant steady state $(u_*, v_*)$  is asymptotically stable for $\tau<\tau^*(d_{21})$ and unstable for $\tau>\tau^*(d_{21})$.
 Fig.\ref{Fig1}  illustrates the stability region and the Hopf bifurcation curves in the $d_{21}-\tau$ plane for $3\leq d_{21}\leq 8$ and $0\leq\tau\leq 6$.  The Hopf bifurcation curves $\tau=\tau_{1,0}$  and $\tau=\tau_{2,0}$ intersect at the point $P_1(4.1354,4.0292)$, which is the double Hopf bifurcation point.  For the point $P_1(4, 2)$ located in the stability region, Fig.\ref{Fig2} illustrates the evolution of the solution of system \eqref{Exa-PP} starting from the initial values  $u_0(x)=0.2+0.1\cos(x/2),~v_0(x)=2.5+0.1\cos(x/2)$, finally converging to the constant steady state $(u_*, v_*)$.

 From Fig.\ref{Fig1}, it is obvious to see that
 $$
\tau^*(d_{21})=
\left\{
\begin{array}{ll}
\tau_{1,0}, & d_{21}^{(1)}<d_{21}<4.1354,\vspace{0.2cm}\\
\tau_{2,0}, &  4.1354<d_{21}<8.
\end{array}
\right.
 $$

\begin{figure}
\centering
\subfloat[]{\includegraphics[scale=0.3]{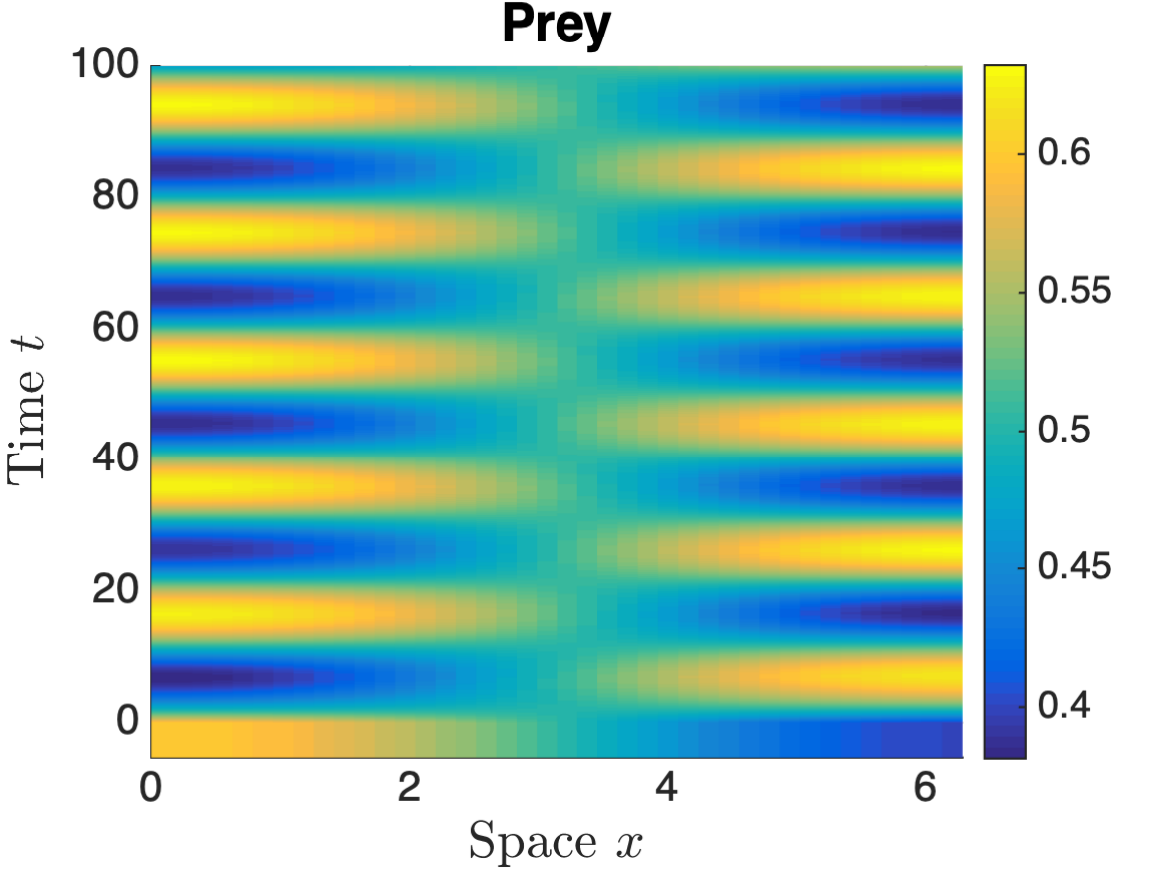}}
\subfloat[]{\includegraphics[scale=0.3]{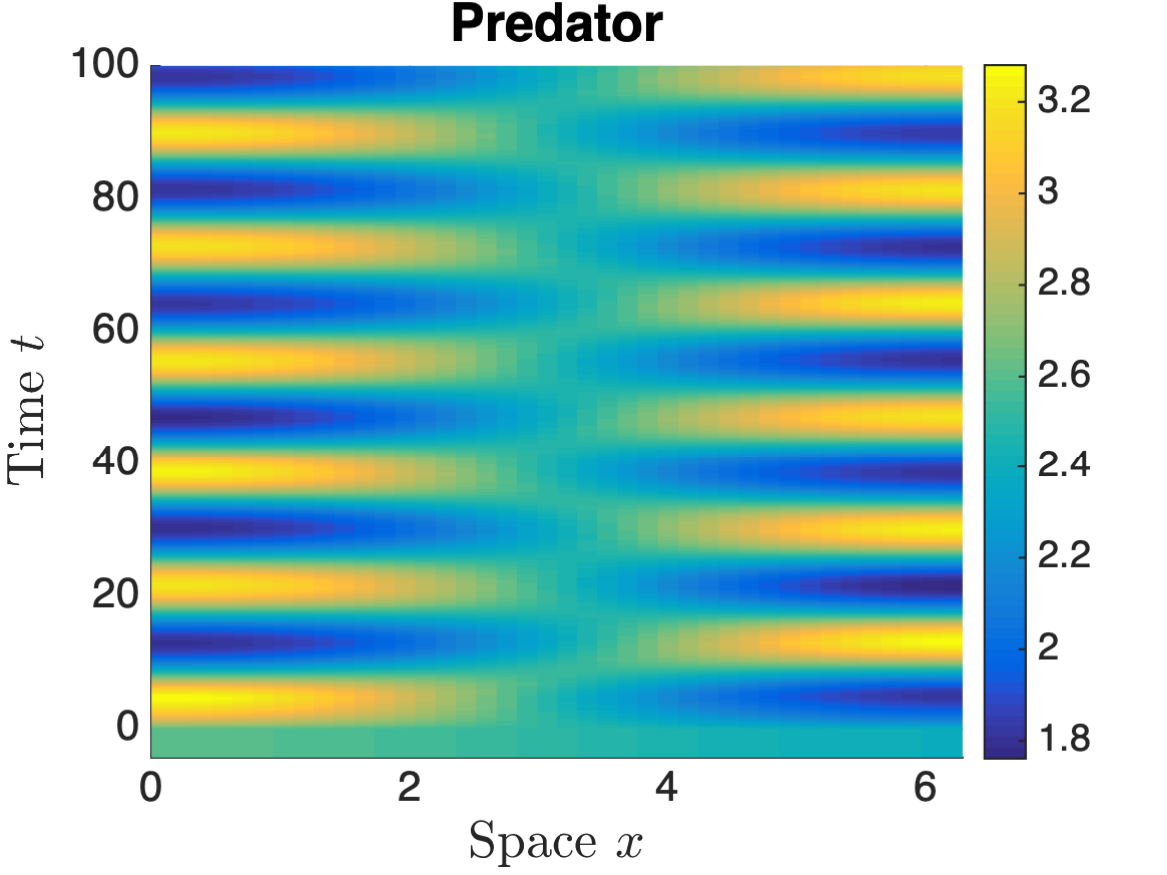}}\quad
\subfloat[]{\includegraphics[scale=0.3]{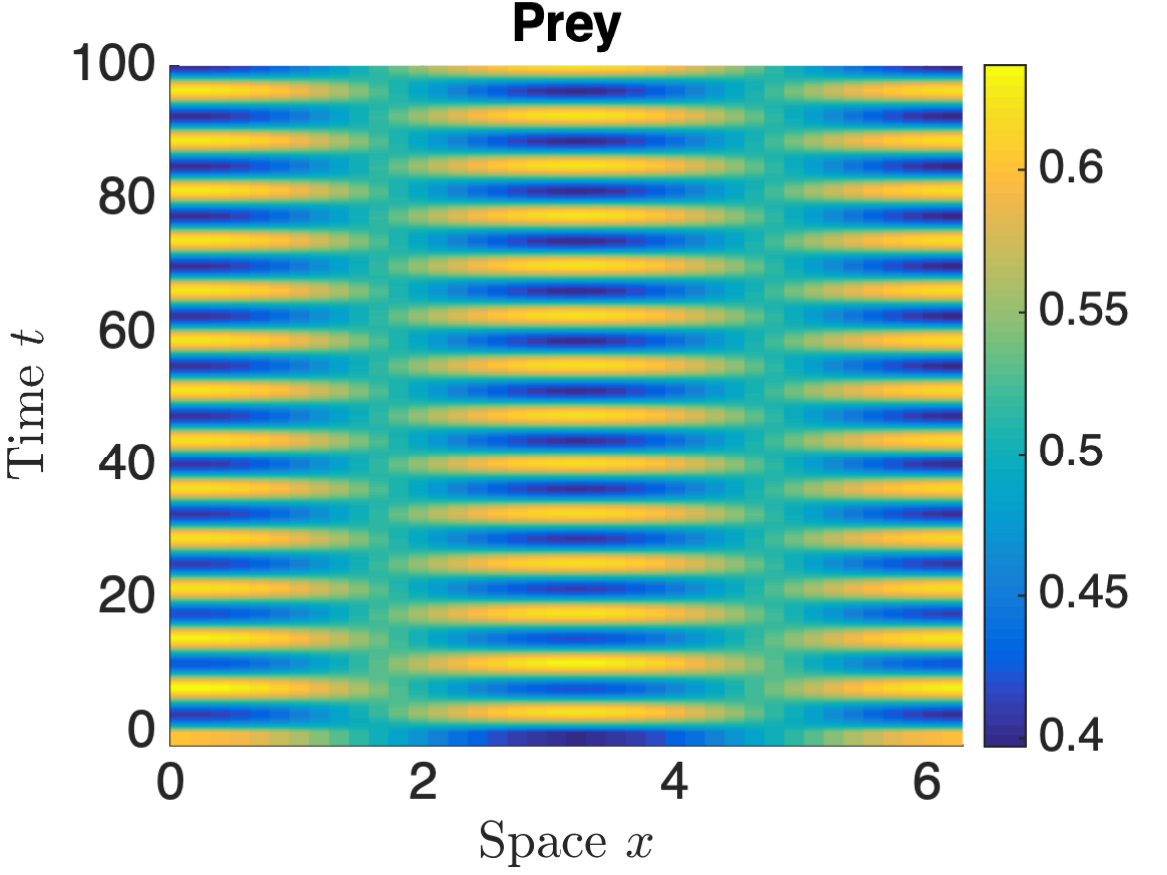}}
\subfloat[]{\includegraphics[scale=0.3]{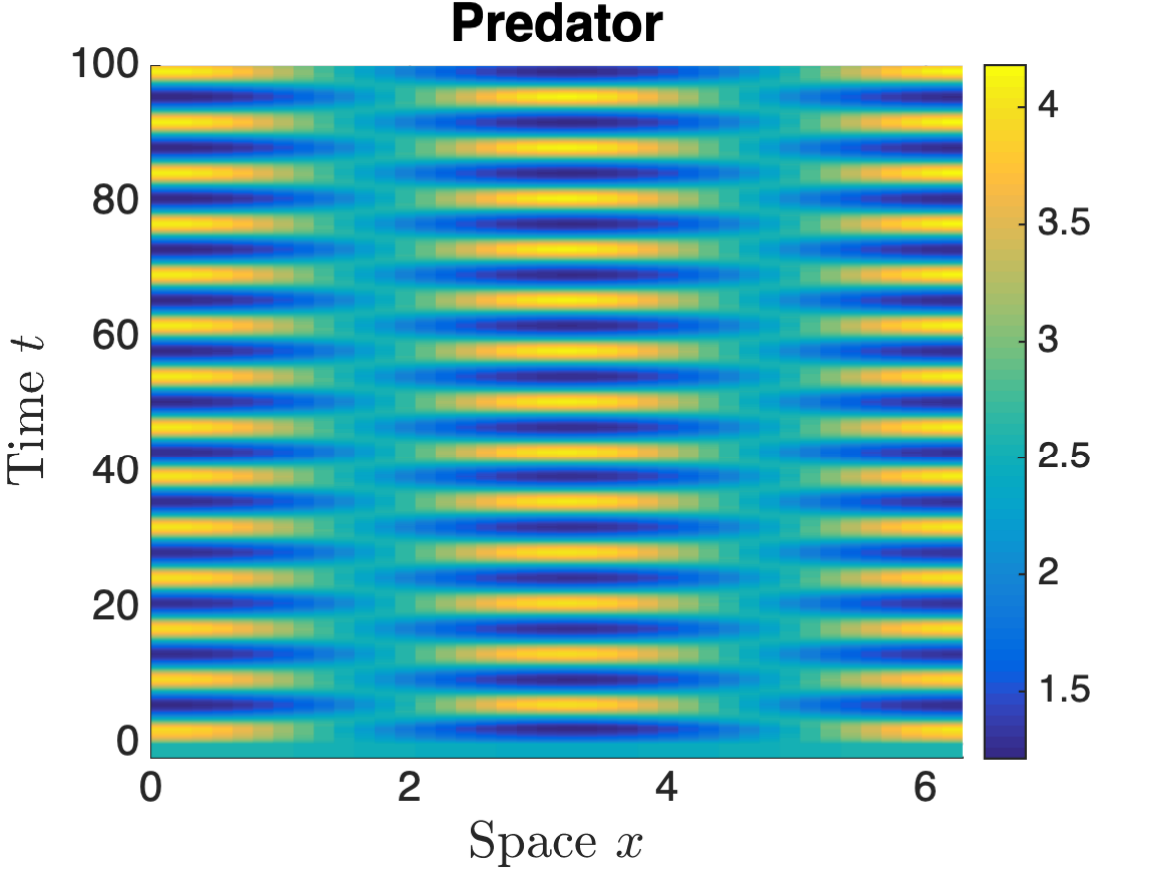}}
  \caption{The spatial-temporal dynamics of system \eqref{Exa-PP}  with the parameters $a=1, ~b=\frac{3}{10}, ~c=\frac{1}{10}, ~d_{11}=\frac{6}{10}, ~d_{22}=\frac{8}{10},~\ell=2$.  $(a)-(b)$ For the  point $P_2(3.6, 5.3)$ of Fig.\ref{Fig1},  there exists  a spatially inhomogeneous periodic solution with mode-1 spatial pattern. The initial values are $u_0(x)=0.2+0.1\cos(x/2),~v_0(x)=2.5+0.1\cos(x/2)$.
  $(c)-(d)$  For the  point $P_3(6, 2)$ of Fig.\ref{Fig1},  there exists  a spatially inhomogeneous periodic solution with mode-2 spatial pattern.  The initial values are $u_0(x)=0.2+0.1\cos(x),~v_0(x)=2.5+0.1\cos(x)$.}
  \label{Fig3}
  \end{figure}

 For   $d_{21}=3.6$, it follows from  \eqref {OMGN}  and   \eqref{TAUNJ} that
 $$\tau_{1,0}\doteq 5.1033 <\tau_{2,0}\doteq 6.6493.$$
For the critical mode-$1$ Hopf bifurcation at $\tau_{1,0}\doteq 5.1033$,  the direction and stability of this Hopf bifurcation  can be determined by the procedure in the previous section with $\tau_c=\tau_{1,0}\doteq 5.1033$ and $d_{21}^c=3.6$. A directional calculation shows that
 $$K_1\doteq 0.0597>0, ~~K_2\doteq -1.5624<0,$$
 which implies that the  mode$-1$  spatially imhomogeneous Hopf bifurcation at  $\tau_{1,0}$ is supercritical and  stable. For $\tau=5.2>\tau_{1,0}$,  Fig.\ref{Fig3}$(a)-(b)$ illustrate the existence of the  spatially inhomogeneous periodic solution with mode-$1$ spatial pattern.

 For  $d_{21}=6$, it follows from  \eqref {OMGN}  and    \eqref{TAUNJ} that
 $$\tau_{2,0}\doteq 1.8398 <\tau_{1,0}\doteq 2.3542.$$
For the critical Hopf bifurcation value $\tau_{2,0}\doteq 1.8398$,   we have
 $$K_1\doteq 0.1733>0, ~~K_2\doteq -2.2283<0,$$
 which implies that the  mode$-2$  spatially imhomogeneous Hopf bifurcation at  $\tau_{2,0}$ is also supercritical and  stable.  For $\tau=2>\tau_{2,0}$,  Fig.\ref{Fig3}$(c)-(d)$ illustrate the existence of the  spatially inhomogeneous periodic solution with mode-$2$ spatial pattern.

\begin{figure}
\centering
 {\includegraphics[scale=0.4]{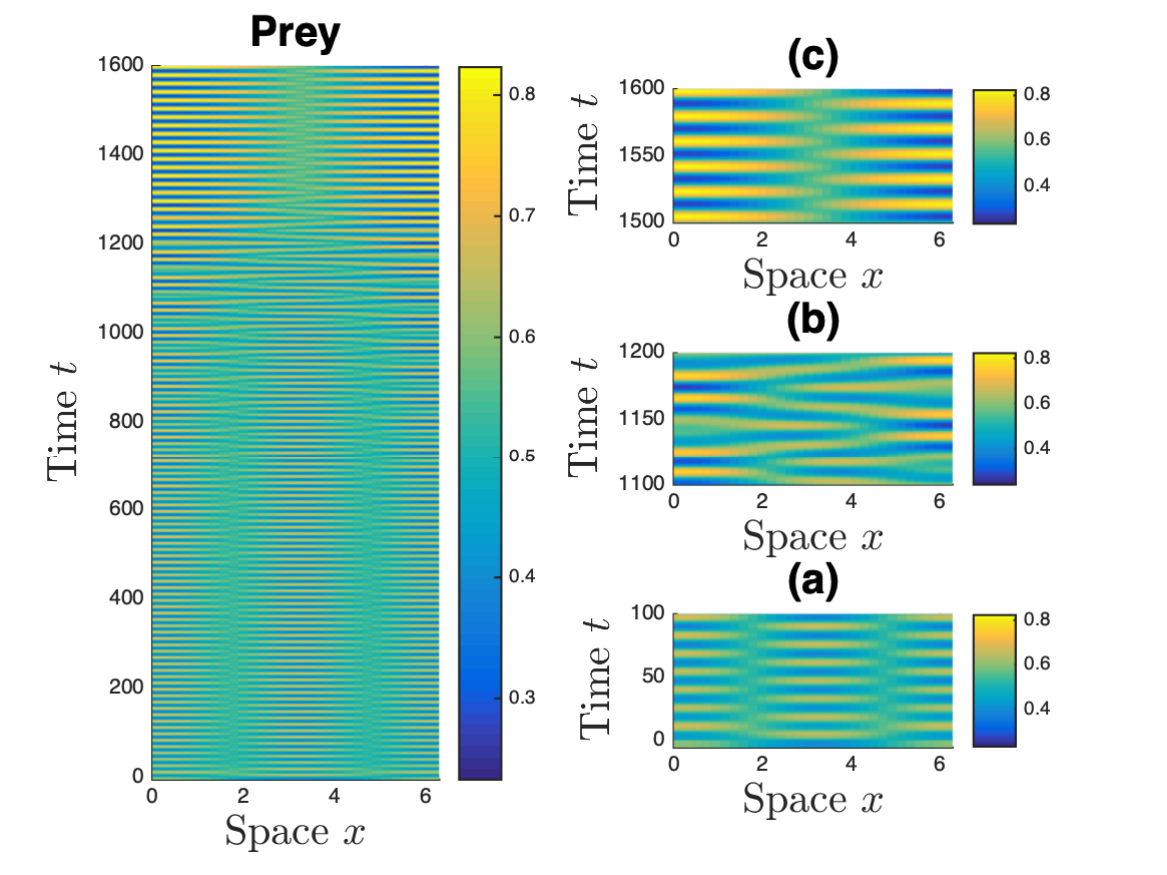}
 \includegraphics[scale=0.4]{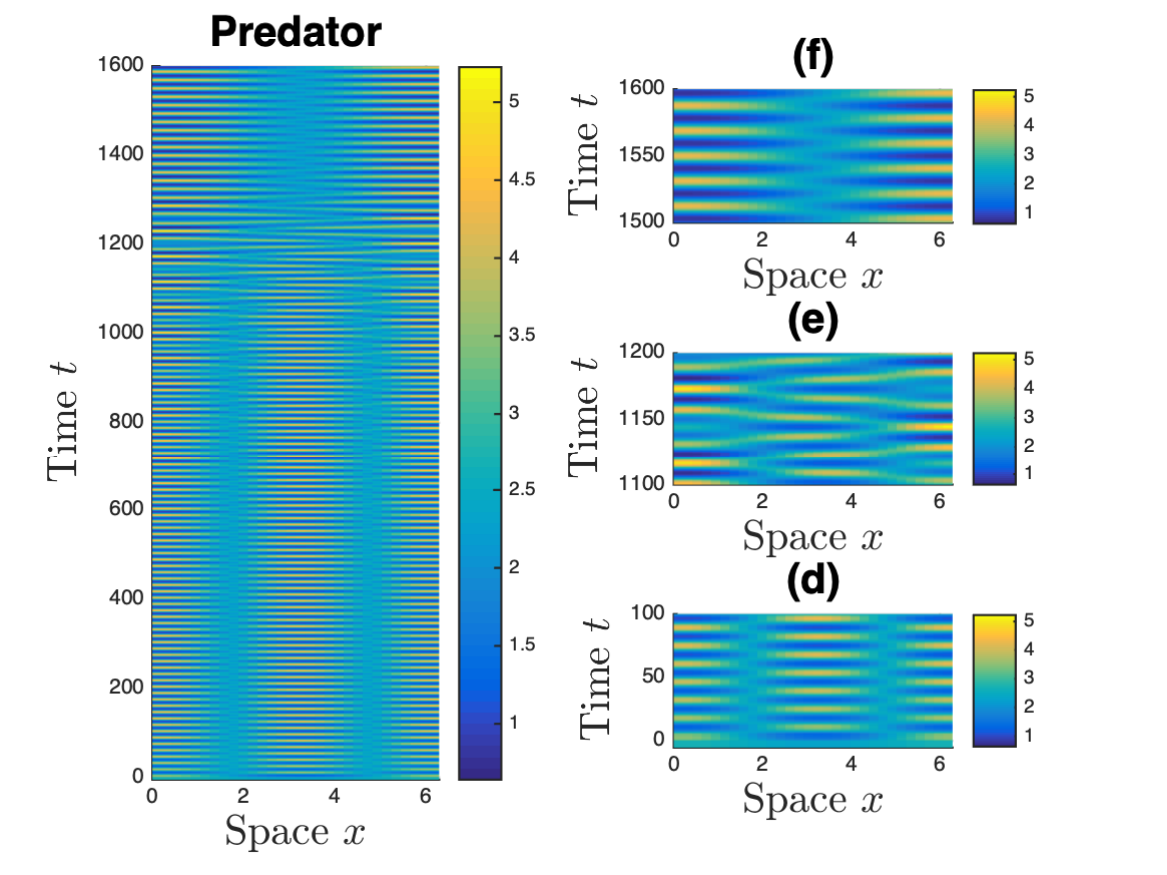}}
  \caption{The spatial-temporal dynamics of system \eqref{Exa-PP}  with the parameters $a=1, ~b=\frac{3}{10}, ~c=\frac{1}{10}, ~d_{11}=\frac{6}{10}, ~d_{22}=\frac{8}{10},~\ell=2$ and  $\left(d_{21}, \tau\right)$ being the point $P_4(4.3, 5.2)$ of  Fig.\ref{Fig1} far from the Hopf bifurcation curves. Pattern transition from the  mode-$2$ spatially inhomogeneous periodic solution   to mode-$1$ spatially inhomogeneous periodic solution. The initial values are $u_0(x)=0.2+0.1\cos(x),~v_0(x)=2.5+0.2\cos(x)$.}
  \label{Fig4}
  \end{figure}

We would like to mention that  the interaction of  mode-$1$ and mode-$2$ Hopf bifurcations leads to more complex dynamics.    To investigate the dynamical classification near the double Hopf point $P$ in detail, the normal form for double Hopf bifurcation should be further developed.   For the point $P_4(4.3, 5.2)$ of  Fig.\ref{Fig1} far from the Hopf bifurcation curves $\tau=\tau_{k, 0}, k=1, 2$, Fig.\ref{Fig4} shows the pattern transition from the  mode-$2$ spatially inhomogeneous periodic solution (Figs.\ref{Fig4}$(a)$ and $(d)$  )  to mode-$1$ spatially inhomogeneous periodic solution (Figs.\ref{Fig4}$(c)$ and $(f)$). Figs.\ref{Fig4}$(b)$ and $(e)$ also illustrate the transiently unstable quasi-periodic patterns.

\section{Conclusion and discussion}\label{sec4}

In this paper, we have developed an algorithm for computing the normal form of the Hopf bifurcation for the reaction-diffusion systems with memory diffusion. Because of  the nonlinearity of the diffusion term  and the the presence of the delay in the diffusion term, the traditional  algorithms for computing the normal form of the Hopf bifurcation for the reaction-diffusion system is not applied to this system. To fill this gap, we generalized the existing algorithms for the reaction-diffusion systems where the diffusion terms are linear and the delay only appears in the reaction terms, to the case where the diffusion terms are not linear  and the delay appears not only in the reaction terms but also in the diffusion terms.

As an illustration of this newly developed algorithm,  we considered a  diffusive predator-prey system with memory-based diffusion and Holling type-II functional response. The memory delay-induced spatially inhomogeneous Hopf bifurcations and the double Hopf bifurcations due to their interactions  are observed.

We  would also like to mention that   the delay-induced double Hopf bifurcation often occurs  in reaction-diffusion systems with memory diffusion and may lead to more complex dynamics like two/three invariant torus.   To determine the dynamical classifications near the double Hopf bifurcation points, new algorithms for computing the normal form of the double Hopf bifurcation of the reaction-diffusion systems with memory diffusion are desired.


\end{document}